\documentclass{amsart} 
\usepackage{amsfonts}
\usepackage{amsmath}
\usepackage{amssymb}
\usepackage{latexsym}
\usepackage{graphicx}
\usepackage{color}
\usepackage{fancyhdr}

\newtheorem{Pro}{Proposition}[section]
\newtheorem{Th}[Pro]{Theorem}
\newtheorem{Cor}[Pro]{Corollary}
\newtheorem{Lem}[Pro]{Lemma}

\theoremstyle{definition}
\newtheorem{Def}[Pro]{Definition}

\theoremstyle{remark}
\newtheorem{Rq}[Pro]{Remark}
\newtheorem{remarks}[Pro]{Remarks}

\newcommand{\field}[1]{\mathbb{#1}}
\newcommand{\Q}{\field{Q}}

\newcommand{\Z }{\field{Z}}
\newcommand{\N }{\field{N}}
\newcommand{\alg}{{\rm{alg}}}
\newcommand{\Gal}{{\rm{Gal}}}

\newcommand{\ov}{\overline}
\newcommand{\wi}{\widetilde}
\newcommand{\gr}{{\rm{gr}}}
\newcommand{\Cq}{q}
\newcommand{\la}{\lambda}

\begin {document}

\title
[Subfields of Nondegenerate Tame Semiramified
Division Algebras]
{Subfields of Nondegenerate Tame Semiramified 
\\ 
Division Algebras}
\author{Karim Mounirh and A. R. Wadsworth}
\address{ }
\email{akamounirh@hotmail.com}
\address{{\rm(for A. R. Wadsworth)} \ Department of Mathematics\\
University of California, San Diego\\
La Jolla, CA   \  92093-0112\\
USA}
\email{arwadsworth@ucsd.edu}
%\date {\today}

\maketitle
\begin{abstract}
  {\normalsize } We show in this article that in many cases the subfields of
  a nondegenerate tame semiramified division algebra of prime power degree over a
  Henselian valued field are  inertial field extensions of
  the center [Th.~\ref{normalsubfield}, Th.~ \ref{pexp} and Prop.~
  \ref{maxcyclic} ].   \\\\
\end {abstract}
\large

{\centerline { \textsc{Introduction}}} 
\bigskip

In their work on abelian crossed product algebras \cite{AS78},
Amitsur and Saltman 
defined a condition they called {\it nondegeneracy} for matrices encoding
the multiplicative structure of such algebras.  They used 
nondegenerate generic abelian crossed product algebras to prove 
the existence of noncyclic $p$-algebras, thereby settling a question 
that had been open since Albert's work on $p$-algebras in the 
1930's. Saltman at that time also showed in \cite{S78} that nondegenerate
generic abelian crossed product $p$-algebras had only one Galois group 
occurring for maximal subgroups Galois over the center, and he
used this to prove the existence of noncrossed product $p$-algebras.
Later, in \cite[Th.~7.17]{S99} he used nondegenerate generic 
abelian crossed products to give examples of indecomposable 
division algebras of exponent $p$ and degree $p^2$,  for any odd
prime $p$, over a field containing a primitive $p$-th root of unity.  
More recently, McKinnie in \cite[Def.~2.4]{Mc07} defined 
nondegeneracy for inertially split semiramified divsion algebras 
over Henselian valued fields  in terms of nondegeneracy of certain 
matrices over the residue 
field; she used this to study prime-to-$p$ extensions of generic 
crossed product $p$-algebras.  In  \cite{Mc08} she further proved 
the indecomposability of nondegenerate inertially split semiramified 
$p$-algebras over Henselian fields of characteristic $p$.
Independently of McKinnie's work, the first author defined 
in  \cite{M08}
nondegeneracy in the somewhat more general context of 
inertially split division algebras over Henselian fields;
in the semiramified case considered in  \cite{Mc08} this 
definition agrees with McKinnie's definition, and also that 
of Amitsur and Saltman.  He proved in particular 
in \cite[Th.~3.5]{M08} that in all 
characteristics a 
nondegenerate inertially split semiramified division algebra of 
prime power degree over a Henselian valued field is indecomposable.

The various formulations of nondegeneracy will be reviewed 
at the beginning of \S 2 below.

In all the work just described, the nondegeneracy condition 
for the algebras was crucial in obtaining constraints on the 
possible subfields of the algebras which are
normal over the center.  Thus, it seems worthwhile to 
investigate more closely the possible subfields of a nondegenerate 
division algebra, particularly the normal subfields.
  We do this here in the Henselian setting,
i.e., where $E$ is a field with a Henselian valuation 
$v$, and $D$ is a division algebra of prime power degree over $E$;
$v$ extends uniquely to a valuation $w$ on $D$, and it is assumed
that $D$ is inertially split and semiramified over $E$ with 
respect to~$w$.  (The valuation-theoretic terminology used here
will be recalled later in this Introduction.)  There is then a 
unique up to isomorphism maximal subfield $N$ of $D$ which is 
inertial (\,= unramified) over~$E$, and $N$~is abelian 
Galois over $E$.  The inertial field extensions of $E$ in~ 
$D$ are fully classified up to isomorphism  as the subfields of 
$N$, and they are all abelian Galois over~$E$.  The question
is thus what other subfields of~$D$ may exist.  This will be 
studied in \S2 below. 

Whenever there is a valuation $v$ on a division algebra $D$, the 
filtration of  $D$ induced by $v$ yields an associated graded ring 
$GD$ in which every nonzero homogeneous element is a unit---this
is called a graded division ring.  When the valuation on the center 
$Z(D)$ is Henselian, the structure of the graded ring closely mimics 
that of $D$.  Likewise, algebraic field extensions of a 
Henselian valued field correspond to graded field extensions of 
the associated graded field.  In \S 1 we will prove some properties 
for algebraic extensions of valued and graded fields, which have some
interest in their own right  and are needed for~\S 2.  We show 
in Th.~\ref{GMnormal} that if $(E, v)$ is a valued field
and $(M, w)$ is a normal finite-dimensional valued
 field extension of $(E, v)$, then 
the associated graded field $GM$  is a normal graded 
field extension of $GE$.
We give in Th.~\ref{lambda} an extension of Hensel's Lemma for 
polynomials over a  valued field all of whose roots have the same
value. We   prove also in Prop.~\ref{wild} that if $E$ is a Henselian valued field with
residue characteristic $p>0$ and $L$ is a  purely wild (resp.,
simple purely wild) finite-dimensional graded field extension of
$GE$, then there is a defectless field extension (resp., a
defectless simple field  extension) $K$ of $E$ such that $GK=L$.
Moreover, if  ${\rm char}(E)=p$, then  $K$ can be  a  purely
inseparable field extension of $E$. We  give in Cor.~ \ref{corresp}
a new (and more explicit) proof of \cite[Th.~5.2]{HW(a)99} which for
an arbitrary Henselian valued field $E$ establishes a  one-to-one
correspondence between the set of isomorphism classes of
finite-dimensional tame field extensions of $E$ and the set of
isomorphism classes of finite-dimensional tame graded field
extensions of $GE$.   We prove also in the last
part of \S 1  some results concerning cyclic graded field
extensions.

In \S 2 we consider subfields of nondegenerate algebras.
Let $E$ be a field with a Henselian valuation ~
$v$, with residue field $\ov E$ and value group 
$\Gamma_E$.  Let $D$ be a 
nondegenerate inertially split 
 semiramified division algebra 
with center $E$, of degree $p^n$ for some 
prime $p$.
 We show  in Th.~ \ref{normalsubfield} that if 
$\text{char}(\ov E) = p$ and  $\Gamma_D/\Gamma_E$~ is not cyclic, 
then any
 subfield of $D$ normal over $E$
is an abelian  Galois inertial extension of~
$E$, 
 and that all  maximal subfields of $D$ Galois over $E$ have the
same Galois group. We prove also  in Prop.~ \ref{elabelian}
that for  $E$ of any  residue characteristic and 
$\Gamma_D/\Gamma_E$  noncyclic,  $D$ is an elementary
abelian crossed product if and only if $\Gal(\overline D/\overline
E)$ is elementary abelian.
Further, we prove in Prop.~ \ref{rk3}
that if  ${\rm rk}(\Gamma_D/\Gamma_E)\geq
3$, then all  subfields of $D$ 
abelian Galois over $E$ are inertial over $E$.
 We show also that if
 ${\rm rk}(\Gamma_D/\Gamma_E)$ is arbitrary but 
${\rm exp}(\Gamma_D/\Gamma_E)=p$,
then any   non-maximal subfield of $D$ is inertial over~$E$. 
In this case, we show in Th.~ \ref{pexp} that if
  $\Gamma_D/\Gamma_E$ is noncyclic and  $K$  is a 
 maximal
  subfield of~ $D$ which is 
normal over $E$ with Galois group not the quaternion group,
then either  $K$  is  cyclic Galois over
$E$ with $[K:E] = p^2$
  or $K$ is inertial over
  $E$. More results concerning subfields  in the case where
   ${\rm char}(\overline E)\nmid{\rm deg}(D)$ are proved at the 
end of \S 2.

\smallskip 
We now
recall some basic  terminology from the theory of 
valued and graded division algebras which will be used 
throughout the paper.

Let  $E$ be a field, $D$ a finite-dimensional division algebra 
over $E$, and $\Gamma$  a
totally ordered abelian group. Let $\infty$ be an element of a set
strictly containing $\Gamma$ with $\infty\notin \Gamma$; extend 
the order on $\Gamma$ to
$\Gamma\cup \{\infty\}$ by setting $\gamma<\infty$ for all
$\gamma\in \Gamma$, and define $\gamma +\infty=\infty+\infty=\infty$.
A map $v\colon D \rightarrow \Gamma\cup \{\infty\}$ is called a
valuation on $D$ if it satisfies the following conditions (for all
$c, d\in D$):
\begin{enumerate}
\item[(1)] 
$v(c)=\infty$ if and only if $c=0$;
\item[(2)] 
$v(cd)=v(c)+v(d)$;
\item[(3)]
$v(c+d)\geq \min\{v(c), v(d)\}$.
\end{enumerate}
We will say that $(D, v)$ is a valued division algebra over $E$.
The value group of $v$ on $D$ is $\Gamma_D = v(D^*)$, where 
$D^* = D\setminus\{0\}$, the group of multiplicative units of $D$.
The residue division algebra is 
$$
\overline D \ = \ V_D\big/ M_D \ = \  \{d\in D \,|\ v(d) \ge 0\}\,
\big/ \,\{d\in D \,|\ v(d) > 0\}.
$$
Of course, $v$ restricts to a valuation on $E$; we write 
$|\Gamma_D:\Gamma_E|$ for the ramification index of 
$D$ over $E$, which is the index in $\Gamma_D$ of its subgroup 
$\Gamma_E$.  Also, we identify the residue field $\ov E$ with its canonical 
image in $\ov D$, and write $[\ov D: \ov E]$ for the residue 
degree of $D$ over $E$, which is the dimension of $\ov D$
as an $\ov E$-vector space.  For background on valued 
division algebras, the reader can consult \cite{JW90} or 
\cite{W}.

Recall the Fundamental Inequality: 
$$
[\overline D : \overline E] \,\, |\Gamma_D : \Gamma_E| \ \le
 \ [D : E] \ < \ \infty.
$$
 $D$ is said to be 
{\it defectless} over $E$ if ${[D : E]=[\overline D : \overline
E]\,|\Gamma_D : \Gamma_E|}$. We say $D$  
{\it inertial} (or {\it unramified})
over $E$ if $[D : E]=[\overline D :
\overline E]$ and the center $Z(\ov D)$ is separable over~ 
$\ov E$. At the other 
extreme, $D$ is  {\it totally ramified} over~$E$ if 
$[D : E]=|\Gamma_D : \Gamma_E|$. 
Now assume $E =
Z(D)$. Then,  $D$~is  {\it inertially split}  if 
it has a maximal subfield which is inertial over $E$. 
Also,  $D$ is said to be  {\it semiramified}  if it is
defectless over~$E$, $\overline D$ is a field,  and $[\overline D :
\overline E]=|\Gamma_D : \Gamma_E|$. It is called 
{\it nicely
semiramified}  if it is inertially split  and 
it has a {\it totally ramified of radical
type}  maximal subfield, i.e., a  maximal subfield
$K$ totally ramified over $E$
such that $K=E[t_1^{1/n_1},..., t_r^{1/n_r}]$, where $t_1,...,
t_r \in E^*$, $\Gamma_K/\Gamma_E=\bigoplus_{i=1}^r\langle
v(t_i^{1/n_i})+\Gamma_E\rangle$ and ${\rm
ord}(v(t_i^{1/n_i})+\Gamma_E)=n_i$. 
%%  If in addition $E$ is
%%  Henselian, then $D$ is said to be  {\it degenerate} if there is an inertial
%%  field extension $M$ of $E$ in $D$ such that the centralizer 
%%  algebra $C_D(M)$ is nicely semiramified and the rank of the
%%  finite abelian Galois group $\Gal(\overline D/\overline M)$ is greater
%%  than 2. Otherwise, we say that $D$ is {\it nondegenerate}.
For $a\in V_D$, we write $\ov a$ for the image of~$a$ in 
$\ov D = V_D/M_D$.  
There is a well-defined canonical group homomorphism
$\theta_D$ mapping   $\Gamma_D/\Gamma_E$ to the Galois group 
$\Gal(Z(\overline D)/\overline E)$, given  by  
$\theta_D(\gamma+\Gamma_E)\colon \ov a \mapsto
\overline {dad^{-1}}$ for all $a\in V_D$ with $\ov a \in Z(\ov D)$, 
where 
 $d$ is  an arbitrary element of $D^*$
with  $v(d)=\gamma$. By \cite[Prop.~ 1.7, Lemma~5.1]{JW90},
$Z(\ov D)$ is normal over $\ov E$ and  
 $\theta_D$ is   surjective; moreover, if $D$ is inertially 
split, then $Z(\ov D)$ is Galois over $\ov E$ and 
$\theta_D$ is an isomorphism.   If  the valuation on $E$ is 
Henselian,   we say
that $D$ is {\it tame}  if 
$D$ is defectless over $E$,  $Z(\overline D)$  is separable  over
$\overline E$, and  $\text{char}(\overline E)
\nmid|\ker(\theta_D)|$.  Equivalently
(see \cite[Prop.~4.3]{HW(b)99}), $D$ is tame iff it is split by the maximal 
tamely ramified field extension of $E$.

We will be working with graded division rings and fields as well
as valued ones.  We recall some of the  terminology and basic
facts in the graded setting, and the connections between the
valued setting and the graded setting.

Let $F$ be an associative ring (with $1$) and 
let $\Gamma$ be a
totally ordered abelian group. We say  that $F$~is   a 
{\it graded ring
of type $\Gamma$} if there are additive subgroups $F_{\gamma}$ ($\gamma \in
\Gamma$) of $F$ such that $F=\bigoplus_{\gamma \in \Gamma}F_{\gamma}$
and $F_{\gamma}F_{\delta} \subseteq F_{\gamma +\delta}$, for all
$\gamma,$ $ \delta \in \Gamma$. In this case, the set
$\Gamma_F=\{\gamma \in \Gamma$ $|$ $ F_{\gamma} \neq 0 \}$ is
called the {\it support}  of~$F$.
If $F$ is a graded ring of type $\Gamma$ and $x \in F_{\gamma}$
for some $\gamma\in \Gamma_F$, we say that $x$ is a homogeneous
element of~$F$; if  $x\ne 0$,  we say that $x$ has {\it grade} 
$\gamma$ and we write
$\text{gr}(x)=\gamma$. 

A graded ring $F$ (of type $\Gamma$) which is
commutative and for which all nonzero homogeneous elements are
invertible is called a {\it graded field}.  Note that, because of
the total ordering on $\Gamma$, in a graded field $F$ 
 every element of the group $F^*$ of multiplicative units 
must be homogeneous; so $F^*$ coincides with the set of
nonzero homogeneous elements of $F$.  Furthermore, the total 
ordering on $\Gamma$ implies that $F$  
 is an integral domain.  It is easy to see also that 
if $M$ is a graded $F$-module (i.e., $M =\bigoplus_{\gamma
\in \Gamma} M_\gamma$ with $F_\gamma M_\delta \subseteq 
M_{\gamma +\delta}$ for all $\gamma,\delta\in \Gamma$)
then $M$ is a free $F$-module with a homogeneous base, and
any two such bases have the same cardinality.  We therefore
write $\dim_F(M)$ for the rank of $M$ as  free $F$-module. 

Let $F$ be a commutative graded ring of type $\Gamma$. An
algebra $A$ over $F$  is called a {\it graded  algebra} (of type
$\Gamma$) over $F$ if $A$ is  a graded ring of type $\Gamma$ and
$F_{\gamma}\subseteq A_{\gamma}$, for all $\gamma \in \Gamma$. 
If $F$ and $A$ are graded fields, we call $A$ a graded field 
extension 
of $F$. If $F$~is a graded field, then a graded algebra over
$F$  in which every nonzero
homogeneous elementis a unit is called a
{\it graded division
algebra}  over $F$. If $F$
is the center of a graded division algebra $A$, then $A$ is called
a {\it graded central division algebra}  over $F$.  We write 
$[A:F]$ for $\dim_F(A)$. For a graded division algebras $A$,
the support set $\Gamma_A$ is a  subgroup of $\Gamma$,
and $A_0$ is a division ring which is an algebra over $F_0$.
Furthermore, it is easy to prove the Fundamental Equality:
$$
[A:F] \ =\  [A_0:F_0]\,|\Gamma_A:\Gamma_F|. 
$$

Let $F$ be  a graded field, let $q(F)$ be its quotient field, and
let $q(F)_{\alg}$ an algebraic closure of $q(F)$.   
Clearly,
for any element $\lambda$ of the divisible hull $\Delta_F$   of
$\Gamma_F$ (i.e., $\Delta_F = \Gamma_F\otimes_{\Z}\Q$\,), there is 
a unique
grading of type $\Delta_F$ on the polynomial ring $F[X]$ extending 
the grading of $F$
and for which $X$ is a homogeneous element with ${\rm
gr}(X)=\lambda$. We denote $F[X]$ with this grading by 
$F[X]^{(\lambda)}$. A
polynomial $f \in F[X]$  is called $\lambda$-{\it homogenizable} 
if 
%% there is  $\lambda \in \Delta_F$ such that 
$f$ is homogeneous in 
$F[X]^{(\lambda)}$.
Let $x \in q(F)_{\alg}$ and let $f_{x, q(F)}$  denote its
minimal polynomial over $q(F)$. We say that $x$ is 
{\it gr-algebraic}
over~$F$ if $f_{x, q(F)}$~is a homogenizable polynomial of
$F[X]$. 
It is shown in \cite[Prop.~2.2]{HW(a)99} that $x$ is gr-algebraic 
over $F$ if and only if the ring $F[x]$ is a graded field 
extension of $F$ and $x$ is homogeneous in $F[x]$.
If $K$~is a graded field extension of $F$, we say that $K$
is {\it gr-algebraic} over $F$ if every homogeneous element of $K$ 
is gr-algebraic over $F$. This holds, in particular, whenever
$[K:F] <\infty$, by \cite[Prop.~2.2]{HW(a)99}. 
Let $F_{\text{gr-alg}}=F[\{x \in q(F)_{\alg}\,|\ \text{$x$
is gr-algebraic over $F$} \}]$, then as proved in 
\cite[Cor.~2.7(c)]{HW(a)99}, 
 $F_{\text{gr-alg}}$ is a gr-algebraic graded field extension of 
$F$ which contains every other gr-algebraic graded field extension
of $F$ in $q(F)_{\alg}$.
We call  $F_{\text{gr-alg}}$ \lq the' {\it graded algebraic closure} 
of~$F$.

Let $K$ be a graded field extension of a graded field $F$
with $[K:F]<\infty$.
In analogy with the valuation terminology, 
   $K$ 
is said to be 
{\it totally ramified} over $F$ if $[K : F]=|\Gamma_K : \Gamma_F|$.
We say $K$ is  {\it inertial}  over~$F$ if $[K : F]=[K_0 : F_0]$
and $K_0$ is separable over~$F_0$. 
Also, $K$ is {\it tame} over $F$  
 if $K_0$ is separable
over~$F_0$ and $\Gamma_K/\Gamma_F$ has no $p$-torsion if
$\text{char}(F) = p \ne 0$. 
Further,  $K$ is {\it purely wild} over $F$ if 
$\text{char}(F) = p\ne 0$, $K_0$
is  purely inseparable over $F_0$, and $\Gamma_K/\Gamma_F$ is a
$p$-group. By \cite[Lemma 3.6]{HW(a)99} $K/F$ is purely wild if and
only if $q(K)/q(F)$ is purely inseparable. 
If $A$ is a  graded central division algebra over $F$,
 we say
that $A$ is {\it semiramified} if $A_0$ is
a field and $[A_0 : F_0]=|\Gamma_A : \Gamma_F|$; 
$A$ is {\it nicely semiramified}  if $A$ has
a maximal subfield inertial over $F$ and another
which is 
 totally ramified  over $F$.

If  $F$ is a graded field  and $A$ is a graded division algebra
of type $\Gamma$
finite-dimensional 
over $F$, we denote by $\Cq(A)$ the  algebra of central quotients
of $A$.  So, $\Cq(A) \cong A\otimes _F q(F)$, which is a division
ring  over $q(F)$ with $[\Cq(A):q(F)] = [A:F] <\infty$.
  The graded structure on $A$
and the total ordering on $\Gamma$ induce a canonical valuation
$v$ on $\Cq(A)$ as follows (see \cite[\S4]{B(a)98} or
 \cite[\S4]{HW(b)99}):  For nonzero 
$a = \sum_{\gamma\in \Gamma}a_\gamma \in A$
(with each $a_\gamma\in A_\gamma$) set $v(a)$ to be the least~ 
$\gamma$ for which $a_\gamma\ne 0$.  Then for nonzero $a\in A$, 
$b\in F$, define $v(ab^{-1}) = v(a) - v(b)$.  It is easy to 
check that $v$ is well-defined and is a valuation on $\Cq(A)$
with $\Gamma_{\Cq(A)} = \Gamma_A$ and $\ov{\Cq(A)} \cong A_0$.
Note that this canonical valuation depends not only on $\Gamma$ 
as a group, but also on the choice of ordering on~$\Gamma$.
Let $Hq(F)$ denote the Henselization of $q(F)$ with respect to its 
canonical valuation (see \cite[\S5.2]{EP} or \cite[\S16]{E72}),
and set  $H\Cq(A) =\Cq(A)\otimes_{q(F)}Hq(F)$. If $F = Z(A)$ 
(so $q(F) = Z(\Cq(A))$, it is known by Morandi's Henselization 
theorem
\cite[Th.~2]{Mor89} that $H\Cq(A)$ is a division algebra over 
$Hq(F)$.  The Henselian valuation on $Hq(F)$ has a unique 
extension to a valuation on $H\Cq(A)$, and it is known that
$H\Cq(D)$ is a tame central
division algebra over $Hq(F)$ (see \cite[Cor.~4.4]{B95} or
\cite[Th.~ 5.1]{HW(b)99}).  

Going in the other direction, suppose we start with a  field
$E$  with a valuation $v$. Then,
the filtration of $E$ induced by $v$ yields a  canonical graded
field $GE$. Namely, for $\gamma\in \Gamma$ let $E^{\gamma}=\{x\in E$ $|$ $ v(x)\geq
\gamma\}$ and $E^{>\gamma}$ = $\{x\in E$ $|$ $ v(x)> \gamma\}$.
Obviously, $E^{>\gamma}$  is a subgroup of the additive group
$E^{\gamma}$. So, we can define the factor group
$GE_{\gamma}=E^{\gamma}/E^{>\gamma}$. For $x\in E\backslash\{0\}$,
we denote by $\widetilde x$ the element $x+E^{>v(x)}$ of
$GE_{v(x)}$; for $0\in E$, set $\wi 0 = 0$ in $GE$. One can easily see that the additive group
$GE=\bigoplus_{\gamma\in \Gamma}GE_{\gamma}$ with the multiplication
law defined for homogeneous elements by $\widetilde x\widetilde
y=\widetilde{xy}$, is a graded field.
Similarly, if $D$ is a  valued division algebra  
finite dimensional over a
field $E$, then the analogous filtration of $D$ yields a  
graded division algebra 
$GD = \bigoplus_{\gamma\in \Gamma}GD_{\gamma}$ where  
$GD_\gamma = D^\gamma/D^{>\gamma}$ (see 
\cite[\S4]{B(a)98}
or \cite[\S4]{HW(b)99}).  Note that
$GD_0 = \ov D$ and $\Gamma_{GD} = \Gamma_D$. 
It is easy to see
 that if $F$ is a graded field and
$D$ is a graded central division algebra over $F$, then  $D$ is 
canonically isomorphic as a graded ring to  
$G\Cq(D)$, the associated graded ring of the valued 
division algebra $\Cq(D)$, via the mapping
$x = \sum_{\gamma\in \Gamma} x_\gamma \mapsto 
\sum_{\gamma\in \Gamma} \wi{x_\gamma}$.
Likewise, $D \cong GH\Cq(D)$, the associated graded ring of
$H\Cq(D)$.

It is  known that graded central division algebras over a 
graded field $F$ play an analogous
role to central division algebras over a Henselian valued 
field. Indeed, their equivalence
 classes form a {\it graded Brauer group} \(\mbox{GBr}(F)\), and there 
is a group isomorphism
 \(\mbox{GBr}(F)\) $\rightarrow$ \(\mbox{TBr}(Hq(F))\), where 
\(\mbox{TBr}(Hq(F))\) is the
 tame part of the Brauer group 
\(\mbox{Br}(Hq(F))\) \cite[Th.~ 5.1]{HW(b)99}. 
Conversely, for any Henselian valued
  field~$E$, there is a canonical group isomorphism
  \(\mbox{TBr}(E)\) $\rightarrow$ \(\mbox{GBr}(GE)\) 
\cite[Th.~ 5.3]{HW(b)99}.\\

\pagestyle{myheadings}
\markright{{\sc Sec.~\ref{Grad}: Graded and valued field extensions}}
\section{Graded and valued field extensions}
\label{Grad}
\setcounter{Pro}{0}

\begin{Lem}\label{homog}
Let $F$ be a graded field, take any $\lambda$ in the divisible hull 
of $\Gamma_F$, and let $f$ be a nonzero $\lambda$-homogenizable 
polynomial in $F[X]$.  Then,
\begin{enumerate}
\item[(1)]
If $h\in F[X]$ and $h\,|\, f$, then $h$ is 
$\lambda$-homogenizable.
\item[(2)]
For $g\in F[X]$, $f\,|\, g$ in $F[X]$ if and only if $f\,|\, g$ in 
$q(F)[X]$.  
\item[(3)]
$f$ is irreducible in $F[X]$ if and only if $f$ is irreducible in 
$q(F)[X]$.  When this occurs, $f$ is a prime element of $F[X]$. 
\end{enumerate}
Thus, unique factorization holds for $\lambda$-homogenizable
polynomials in $F[X]$.
\end{Lem}

\begin{proof}
(1) This holds because $\Gamma_{F[X]^{(\lambda)}}$ is totally 
ordered.  Therefore, for nonzero $h,k$ in $F[X]$, the 
lowest [resp.~ highest] grade homogeneous component of $hk$ 
is the product 
of the lowest [resp.~highest] grade components of $h$ and $k$.
So, if $hk$ is homogeneous, then $h$ and $k$ must also be 
homogeneous.  

(2)  Write $f = \sum_{i=0}^n a_iX^i$.  Since $f$ is 
$\la$-homogenizable, each nonzero $a_i$ is homogeneous
in $F$, so 
lies in $F^*$.  (2) thus follows by the division algorithm 
for polynomials, since the leading coefficient of $f$ is 
a unit.

(3)  Since the leading coefficient of $f$ lies in $F^*$, 
we may assume that $f$ is monic.  Because the integral domain
$F$ is integrally closed by \cite[Cor.~1.3]{HW(a)99}, if $f$ is 
irreducible in   $F[X]$, then $f$ is irreducible in $q(F)[X]$.
Conversely, if $f$ is irreducible in $q(F)[X]$, then 
$fF[X]$ is a prime ideal of $F[X]$, since (2) shows that
$fF[X] = (f\,q(F)[X]) \cap F[X]$. Hence $f$ is a prime element
of $F[X]$, so it is irreducible in $F[X]$.  

Since nonzero constant $\la$-homogenizable  polynomials are 
units of $F[X]$,
it follows by induction on degree and by (1) and (3) above that 
every $\la$-homogenizable polynomial of positive degree 
is a product of prime $\la$-homogenizable polynomials.  The 
usual argument gives the uniqueness of such a factorization.  
\end{proof}

Let $F$ be a graded field and 
let $L$ be an algebraic
graded field extension of $F$. Then,  we say that $L$~is {\it normal}
over $F$ if every homogenizable irreducible polynomial $g$ of $F[X]$
which has a root in $L$ factors into polynomials of degree one in
$L[X]$. When this occurs, each root $x$ of such a $g$ is homogeneous
in~$L$, since $X-x$ is homogenizable in $L[X]$  by Lemma~\ref{homog}(1).  
Moreover, the minimal polynomial $f_{x,q(F)}$
of $x$ over~$q(F)$ lies in $F[X]$ (as $L$ is integral over $F$, 
which is integrally closed), and $f_{x,q(F)}$ is $\la$-homogenizable, 
where $\la = \gr(x)$. So by Lemma~\ref{homog}, $g = af_{x,q(F)}$
for some $a\in F^*$.  
Thus,  $L$ is normal over $F$ if and only if for any
$x\in L^*$, $f_{x, q(F)}$ factors into  polynomials of degree one
in $L[X]$. 

\begin {Lem}\label{qnormal}
Let $L/F$ be an algebraic graded field extension. Then, $L$ is 
normal over $F$ if
and only if $q(L)$ is a normal field extension of $q(F)$.
\end {Lem}

\begin{proof} Suppose that $L$ is normal over $F$ and consider a
$q(F)$-monomorphism $\sigma$ from $q(L)$ into $q(L)_{\alg}$,  
the algebraic closure of $q(L)$. Let $x\in L^*$ and let $f_{x,
q(F)}$ be its minimal polynomial over $q(F)$. Obviously, we have
 $f_{x, q(F)}(\sigma(x))=0$. It follows by the normality of $L/F$ that
 $\sigma(x)\in L$. Note that we have  $q(L)=L\cdot q(F)$,  so $\sigma(q(L))=q(L)$. Therefore,
 $q(L)$ is a normal field extension of $q(F)$.

Conversely, suppose that $q(L)$ is a normal field extension of
$q(F)$ and let $g$ be a homogenizable irreducible polynomial of
$F[X]$ with a root in $L$. By Lemma~\ref{homog}(3) $g$ remains 
irreducible over
$q(F)$,  so by the normality $g$ splits over $q(L)$. 
Clearly, the roots
of $g$ are integral over $F$, so they all lie  in $L$ as 
 $L$~is integrally closed.
\end{proof}

\begin {Pro}\label{tandn}
 Let $L/F$ be a finite-dimensional graded  field extension. 
Then, the following
 are equivalent:
\begin{enumerate}
\item[(1)]
$L/F$ is tame and normal.
\item[(2)]
 $L$ is  a Galois graded field extension of $F$.
\end{enumerate}
\end {Pro}

\begin{proof} This follows by  \cite[Th.~ 3.11(a),(b), 
 Lemma 1.2]{HW(a)99}. 
\end{proof}

Let $L/F$ be a normal finite-dimensional graded field extension.
The Galois group  of $L$ over $F$ is the group $\Gal(L/F)$
consisting of graded (i.e., grade-preserving)
$F$-isomorphisms of $L$. Let $G = \Gal(L/F)$ and let
$L^G$ denote the set of elments of $L$
  invariant under the action of $G$; so, 
$L^G$ is a graded subfield of~$L$. 
%% Recall from \cite[????]{HW(a)99} that $L$ is 
%% Galois over $F$ if and only if $L^G = F$.   
It was  proved in \cite[p.~26]{B00} that $L$ is
 tame over $L^G$. The following proposition gives a more
 general result.

\begin {Pro}\label{fix}
 Let $L/F$ be a finite-dimensional normal graded field extension
 with Galois group $G$. Then, $L^G$ is purely 
wild over $F$ and $L$
  is Galois over $L^G$. Moreover, if $T$ is the tame closure 
of $F$ in~$L$, then
  $L=T\cdot L^G\cong_gT\otimes_FL^G$.
\end {Pro} 

\begin{proof} Since $L/F$ is normal , by 
Lemma~ \ref{qnormal} $q(L)/q(F)$ is a normal field
 extension.  By \cite[Cor.~ 2.5(d)]{HW(a)99}, every 
$\sigma \in \Gal(q(L)/q(F))$ restricts to a graded 
$GF$-automorphism of $GL$.  Furthermore, the map
$\Gal(q(L)/q(F)) \to \Gal(L/F) =G$ by restriction is 
an isomorphism as $q(L) = L\otimes_Fq(F)$.  Therefore,
we identify $\Gal(q(L)/q(F))$ with $G$. Recall from 
field theory (see, e.g.,  \cite[Prop.~ 3.2, p.~316]{Karp89}) 
that the normality of  $q(L)$  over $q(F)$, implies $q(L)$ is Galois over 
$q(L)^G$, which is purely inseparable over~$q(F)$.  Moreover,
if $S$ is the separable closure of $q(F)$ in $q(L)$, then 
${q(L) = S\cdot q(L)^G \cong S\otimes_{q(F)} q(L)^G}$.
Now, since every $x$ in $q(L)$ is expressible as $ab^{-1}$
with $a\in L$ and $b\in F\setminus \{0\}$, we have ${q(L)^G = 
L^G\cdot q(F) = q(L^G)}$.  Hence, $L$ is Galois over $L^G$
by  \cite[Th.~ 3.11(b)]{HW(a)99} as $q(L)$ is Galois over $q(L^G)$,
and $L^G$ is purely wild over~$F$ by 
\cite[Lemma 3.6]{HW(a)99} as $q(L^G)$ is purely inseparable over 
$q(F)$.   By \cite[(3.8)]{HW(a)99}, $q(T) = S$.  Because 
$T\otimes_F L^G$ is a  torsion-free $F$-module, 
it injects into ${(T\otimes_F  L^G) \otimes_F q(F) \cong
q(T)\otimes _{q(F)} q(L^G)\cong q(L)}$. Let $L' = T\cdot L^G$, 
which is the image of $T\otimes_FL^G$ under its injective
mapping to $q(L)$.  Then, $L'$ is a graded subfield of $L$, 
and the isomorphism  $T\otimes_F L^G\to L'$ respects the gradings.
Also, ${[L:L'] = [q(L):q(L')] = 1}$, as $q(L') = q(T) \cdot q(L^G)
= q(L)$. Thus, $L = L' =T\cdot L^G \cong_g  T\otimes_F L^G$.  
\end{proof}

\begin{Th}\label{GMnormal}
  Let $(E, v)$ be a valued field and $(M, w)$ a   finite-dimensional valued field
extension of $(E, v)$. If $M$ is normal over $E$, then $GM$ is 
normal over $GE$.
\end{Th}

\begin{proof}  Assume first that $w$ is the unique valuation 
of $M$ 
extending $v$ on $E$. Let $x\in M$ and let $f_{x,
E}$ be its minimal polynomial over $E$. Since $M$ is normal over
$E$,  we can write $f_{x, E}=\prod_{i=1}^n(X-x_i)$, where
$x_i=\sigma_i(x)$ for some $\sigma_i\in \Gal(M/E)$. Moreover, 
since
$w$ is the unique extension of
  $v$ to $M$,  $w(x_i)=w(x)$ for 
$1\leq i\leq n$. Therefore, the polynomial 
$\prod_{i=1}^n(X-\widetilde  {x_i})$ lies 
in $GE[X]$, by \cite[ Lemma~2.1, Lemma~2.4]{B(b)98}  
(or by Lemma~\ref{split} below, which shows that $h = 
f^{(\lambda)} \in GE[X]$).
 Hence, the minimal polynomial $f_{\widetilde  x, q(GE)}$ of 
$\widetilde  x$ over $q(GE)$
 splits into polynomials of degree one in $GM[X]$,
showing that $GM$~is normal over $GE$. 

Now, without assuming that $w$ is the unique extension of $v$ to
$M$, let $I=M^{\Gal(M/E)}$ and let $K$~be the decomposition field of
 $w$ over $I$.
 Since $M$ is normal over $K$ and $w$~is the unique extension of 
its restriction $w|_{K}$ to
 $M$,  by the first part of the proof $GM$ is normal over $GK$. 
So by 
Lemma~\ref{qnormal}, $q(GM)$~is normal
 over $q(GK)$. Moreover, since $(K, w|_{K})$ is an immediate field extension
 of $(I, w|_{I})$ by \cite[Cor.~5.3.8(0), pp.~134--135]{EP}, we have 
  $GK=GI$. Note that because $I$~is purely inseparable over~$E$,
we have $GI$ is purely wild over $GE$, so 
   $q(GI)$ is purely inseparable over $q(GE)$.  
 Therefore, $q(GM)$ is
  normal over $q(GE)$, so again by Lemma ~\ref{qnormal}, $GM$ is 
normal over $GE$. 
\end{proof}

In what follows we will consider polynomials over a valued field
$(E, v)$ for which all the roots in an algebraic closure $E_{\alg}$
of $E$ have the same value for any valuation that extends $v$ to
$E_{\alg}$.  The following proposition  generalizes 
\cite[Lemma 2.1]{B(b)98}, which gives (1) $\Leftrightarrow$ (3) 
under the additional assumptions that $v$~is Henselian and 
$f$ is monic. 

  \begin {Pro}\label{samevalue}
 Let $(E, v)$ be a valued field, $E_{\alg}$ an algebraic closure 
of $E$,  and let ${f=\sum_{i=0}^na_iX^{i}\in E[X]}$ 
with $a_0a_n\neq
  0$. Let $\lambda = \frac 1n \big(v(a_0)-v(a_n)\big)$ in the divisible 
hull of  $\Gamma_E$.
Then, the following statements are equivalent:
\begin{enumerate}
\item[(1)] 
For some extension of $v$ to $E_{\alg}$, all the roots of $f$ 
in $E_{\alg}$ have the same value.
\item[(2)] 
For every extension of $v$ to $E_{\alg}$, all the roots of $f$
in $E_{\alg}$ have  value $\lambda$.
\item[(3)] 
  $v(a_i)\ge 
(n-i)\lambda +v(a_n)$ for every $i$, $0\le i\le n$.
\item[(4)] 
Let $w$ be an extension of $v$ to $E_{\alg}$, $c\in E_{\alg}$
with $w(c)=\lambda$, and  let
$h=\frac{1}{a_nc^n}f(cX)$. Then $h$ is a monic polynomial of
$V_{\alg}[X]$, where $V_{\alg}$ is the valuation ring of $w$.
\end{enumerate}
  \end {Pro}

  \begin{proof} $(2)\Rightarrow (1)$ is clear.

 $(1)\Rightarrow (3)$ Let $f = a_n(X-x_1)\ldots(X-x_n)$
in $E_\alg[X]$, and let $s_j$ be the $j$-th symmetric polynomial
in $x_1, \ldots, x_n$ for $1\le j\le n$.  Suppose the $x_i$ all have
the same value for some extension $w$ of $v$ to~$E_\alg$.  
Then $w(x_i) = \lambda$ for all $i$, as $a_0 = 
(-1)^na_nx_1\ldots x_n$.  Since $s_j$ is a sum of products of $j$ of the 
$x_i$, $w(s_j) \ge j\lambda$.  For $0\le i\le n-1$ we have
$a_i = (-1)^{n-i}s_{n-i}a_n$; hence, $v(a_i)  =
w(s_{n-i}) +v(a_n)\ge (n-i) \lambda +v(a_n)$.     

  $(3)\Rightarrow (2)$ Let $w$ be an extension of $v$ to 
$E_{\alg}$, and let $x$ be any root of $f$ in $E_\alg$.  
Since $\sum _{i=0}^n a_ix^i = 0$, in the list of values
$w(a_0), w(a_1x), \ldots, w(a_nx^n)$ the least value must
occur at least twice.  If $w(x) > \lambda$, then 
(3) yields for $i>0$,
$$
w(a_ix^i)  \ > \  v(a_i) + i\lambda  \ \ge  \ (n-i) \lambda + v(a_n) 
+i\lambda =n\lambda +v(a_n)  \ = \ v(a_0).
$$
This is not possible, since then the least value on the list 
would be  $w(a_0)$, occurring only once. Similarly, if 
$w(x) <\lambda$,
then for $i <n$,
\begin{align*}
w(a_nx^n) \ &= \ v(a_n) +i\, w(x) + (n-i) w(x)   \\ 
&< \ v(a_n) +i\,w(x) +(n-i)\lambda  \ \le  \ 
v(a_i) +i\,w(x)  \ = \ w(a_i x^i).
\end{align*}
This is also ruled out, since  the least 
value on the list would be then $w(a_nx^n)$, occurring only
once.  Therefore, $w(x) = \lambda$ for any root $x$ of $f$.

  $(3)\Leftrightarrow (4)$ Clearly, $h$ is a monic 
polynomial.
  Write $h=\sum_{i=0}^nb_iX^i$, where $b_i=a_ia_n^{-1}c^{i-n}$.
 Then, 
  $w(b_i)=v(a_i)-v(a_n)+(i-n)\lambda$.  Hence, 
$w(b_i) \ge 0$ if and only if 
$v(a_i) \ge  (n-i)\lambda +v(a_n)$.
\end{proof}

\begin{Def}If $f = \sum_{i=0}^na_ix^i\in E[X]$ satisfies 
the equivalent conditions of 
Prop.~\ref{samevalue}, then we call $f$ a 
{\it $\lambda$-polynomial}, 
where $\lambda = \frac{1}{n}\big(v(a_0)-v(a_n)\big)$ is the common
value of all the roots of $f$.
  We then write
 $\widetilde f^{(\lambda)}:=
\sum_{i=0}^n\ov a_i^{(\lambda)}X^i\in GE[X]$,
 where $\ov a_i^{(\lambda)}$ is the class of $a_i$ in
 $GE_{(n-i)\lambda+v(a_n)}$. Observe that 
$\ov
 a_i^{(\lambda)}=\widetilde a_i$ if $v(a_i)=(n-i)\lambda+v(a_n)$
 and $\ov a_i^{(\lambda)}=0$ if
 $v(a_i)>(n-i)\lambda+v(a_n)$. Thus, $\widetilde f^{(\lambda)}$~
is a $\lambda$-homogenizable polynomial in $GE[X]$ with 
${\rm {gr}}(\wi f^{(\lambda)}) = v(a_0)$ and $\deg(\wi f^{(\lambda)})
=\deg(f)$.  
\end{Def}

\begin{Lem}\label{split}
Let $(E,v)$ be a valued field, and let 
$f = \sum_{i=0}^n a_iX^i$ be a $\lambda$-polynomial in $E[X]$.
Let $K$ be an algebraic field extension of $E$ over which
$f$ splits, say $f = a_n(X-x_1)\ldots(X-x_n)$ in $K[X]$, 
and let $w$ be any extension of $v$ to $K$.  
Then, $\wi f^{(\lambda)}= \wi a_n(X- \wi x_1)\ldots (X-\wi x_n)$
in $GK[X]$.  
\end{Lem}

\begin{proof}
We have $f = \sum_{i=0}^n a_iX^i = a_n\prod_{i=1}^n(X-x_i)$.  
Let $s_0 = 1$ and for $1\le k\le n$ let $s_k$ be the $k$-th 
symmetric polynomial in $x_1, \ldots, x_n$.  So, $a_i = 
a_n(-1)^{n-i}s_{n-i}$ for $0\le i\le n$.  Let $ g = 
\wi a_n \prod_{i=0}^n (X-\wi x_i) = \sum_{i=0}^n b_iX^i$ in 
$GK[X]$. Let $t_0 = \wi 1$ and let $t_k$ be the $k$-th 
symmetric polynomial in $\wi x_1, \ldots, \wi x_n$ for 
$1\le k\le n$.  So, each $b_i = \wi a_n (\wi{-1})^{n-i}t_{n-i}$. 
Now, each $s_k$ is a sum of monomials of degree $k$ 
in the $x_i$ (so of
value $k\lambda$).  Hence, $w(s_k) \ge k\lambda$.  
We have $(x_{j_1} \ldots x_{j_k})\wi{\phantom x} =
\wi x_{j_1} \ldots \wi x_{j_k}$ for all indices 
$j_1, \ldots, j_k$.  When $w(s_k) = k\lambda$, $\wi s_k$~is 
the sum of the images of its monomials in $GK$, i.e, 
$\wi s_k = t_k \ne 0$ in $GK_{k\lambda}$.  When $w(s_k) >
k\lambda$, the sum of its images in $GK_{k\lambda}$ is $0$, 
i.e., $t_k = 0$.  Now consider $\wi f^{(\lambda)} = \sum 
_{i=0}^n\ov a_i^{(\lambda)}X^i$. 
If $\ov a_i^{(\lambda)} \ne 0$, then $v(a_i) = (n-i)\lambda +v(a_n)$, 
so $v(s_{n-i}) = (n-i)\lambda$, and hence,
$$
\ov a_i^{(\lambda)} \ = \ \wi a_i \ = \ 
\big(a_n(-1)^{n-i}s_{n-i}\big)\wi{\phantom x} \ = \ 
\wi a_n(\wi{-1})^{n-i}\wi s_{n-i} \ =  \ 
\wi a_n(\wi{-1})^{n-i}t_{n-i} \ = \ b_i.
$$ 
On the other hand, if $\ov a_i^{(\lambda)} = 0$, then
$v(a_i) > (n-i)\lambda +v(a_n)$, yielding  
$v(s_{n-i}) > (n-i)\lambda$,
so $t_{n-i} = 0$; hence, $b_i = \wi a_n(\wi{-1})^{n-i}t_{n-i} = 0
= \ov a_i^{(\lambda)}$.  Thus, $\wi f^{(\lambda)} = g$.  
\end{proof}

Let $(E,v)$ be a Henselian valued field.  The next 
theorem generalizes to  arbitrary 
$\lambda$-polynomials over $E$ 
well-known basic properties  for $0$-polynomials,
which are those $f = \sum_{i=1}^n a_iX^i$ in 
$V_E[X]$ with $v(a_n) = v(a_0)=0$ (cf. ~\cite[Th. 4.1.3,
pp.~87--88]{EP}).

\begin {Th}\label{lambda}  
Let $(E, v)$ be a Henselian valued field,  
$f=\sum_{i=0}^na_iX^{i}$ a polynomial
 of $E[X]$ with $a_0a_n\neq 0$ and let 
$f'=\sum_{i=0}^n\widetilde a_iX^{i}\in GE[X]$. Then, 
\begin{enumerate}
\item[(1)] 
If $f$ is a $\lambda$-polynomial and $f=gh$ in $E[X]$, then $g$ and $h$ are
$\lambda$-polynomials and
$\widetilde f^{(\lambda)}=\widetilde g^{(\lambda)}\widetilde h^{(\lambda)}$ in $GE[X]$. So,
 if $\widetilde f^{(\lambda)}$ is irreducible in $GE[X]$, then $f$ is irreducible in $E[X]$.
\item[(2)] 
If $f$ is irreducible in $E[X]$, then $f$ is a
$\lambda$-polynomial for $\lambda=\frac{1}{n}(v(a_0)-v(a_n))$.
Furthermore, $\widetilde f^{(\lambda)}= \widetilde a_nk^s$ for
some irreducible monic homogeneous polynomial  $k$ of
$GE[X]^{(\lambda)}$ and some positive integer $s$.
\item[(3)] 
If $f$ is a $\lambda$-polynomial in $E[X]$ and if $\widetilde
f^{(\lambda)}=\ell m$ in $GE[X]$ with ${\rm{gcd}}(\ell, m)=1$,
then there exist $\lambda$-polynomials $g$, $h$ in $E[X]$ such
that $f=gh$,  $\widetilde g^{(\lambda)}=\ell$, and $\widetilde
h^{(\lambda)}=m$.
\item[(4)]
If $f$ is a $\la$-polynomial and $\wi f^{(\la)}$ has a simple root 
$b$ in $GE$, then $f$ has a simple root $a$ in $E$ with $\wi a = b$.
\item[(5)] 
 Suppose $f'$ is a $\lambda$-homogenizable   polynomial of
$GE[X]$.
Then, $f$ is a $\lambda$-polynomial and
$\wi f^{(\lambda)} = f'$.
\end{enumerate}
\end {Th}

\begin{proof} (1) If $x$ is any root of $g$ or of
$h$ in $E_\alg$, then $x$ is also a root of $f$, so $x$ has 
value
$\lambda$. Thus, $g$~and~$h$ are $\lambda$-polynomials.
Let $K$ be any algebraic extension of $E$ over 
which $f$ splits (so $g$ and $h$ split), and let
$w$ be any extension of $v$ to $K$.  In $K[X]$ write
$g = b\prod_{i=1}^r (X-x_i)$ and $h = c\prod_{i=r+1
}^n(X-x_i)$.  Then the leading coefficient of 
$f$ is $bc$, and $\wi{bc} = \wi b \wi c$ in $GE$.  
By applying Lemma~\ref{split} to $f$, $g$, and $h$, 
we obtain $\wi f^{(\lambda)} = \wi g^{(\lambda)}\wi h^{(\lambda)}$.

(2) Let $N$ be a normal field extension of $E$ that contains the
roots $(x_i)_{i=1}^n$ of $f$ and denote the unique extension of 
the Henselian $v$ to
$N$  by $w$. Since $f$ is irreducible in $E[X]$, for
 any $1\leq i\leq n$, there is an $E$-automorphism $\sigma_i$ 
of $N$ such that
 $\sigma_i(x_1)=x_i$. Because $v$ is Henselian, $w\circ \sigma_i
= w$. 
So, $w(x_i)=w(x_1)$. This shows that $f$ is a
 $\lambda$-polynomial for $\lambda=w(x_1)$. Moreover, 
as $\sigma_i$ preserves $w$, it
  induces a graded $GE$-automorphism $\widetilde \sigma_i$ on 
$GN$ for which
  $\widetilde\sigma_i(\widetilde x_1)=\widetilde x_i$.
This automorphism of course extends to a $q(GE)$-automorphism
of $q(GN)$.  
So,  
the minimal polynomial $k$ of
  $\wi x_1$ over $q(GE)$ is also the minimal polynomial
of $\wi x_i$. Since the monic irreducible factors of $\wi f^{(\la)}$
in $q(GE)[X]$ are the minimal polynomials of its roots
and Lemma~\ref{split} shows that the roots of $\wi f^{(\la)}$
are the $\wi x_i$, we must have $\wi f^{(\la)} =\wi a_n k^s$ in
$q(GE)[X]$.  This $k$ lies in $GE[X]$ as noted above (because
 the graded field $GE$ is integrally closed), and $k$ is 
$\la$-homogenizable by Lemma~\ref{homog}(1) above.

(3) Note that $\ell$ and $m$ are $\la$-homogenizable by 
Lemma~\ref{homog}(1).
Without loss of generality, we may assume that $f$, $\ell$, 
and $m$ are
 monic polynomials. Write ${f=\prod_{i=1}^rp_i^{t_i}}$, where 
$p_i$ are distinct
  monic irreducible polynomials in $E[X]$. By (1) above the $p_i$ 
are also $\lambda$-polynomials and by (2) each 
$\widetilde p_i^{(\lambda)}$ equals $q_i^{s_i}$ for some monic
  irreducible $\la$-homogenizable
polynomial $q_i$ in $GE[X]$, so 
$\widetilde f^{(\lambda)}=\prod_{i=1}^rq_i^{s_it_i}$. 
Since each $q_i$ is a prime element of $GE[X]$ by 
Lemma~\ref{homog}(3), 
    $q_i$ divides $\ell$ or $m$ but not both. Let 
$g$ be the product of those
   $p_i^{t_i}$ with $q_i$ dividing $\ell$, and 
$h$  the product of 
the $p_i^{t_i}$ with $q_i$ dividing
   $m$. Then, $f=gh$.  Furthermore, 
$\widetilde g^{(\lambda)}=\ell$,  and 
$\widetilde h^{(\lambda)}=m$ by the unique 
factorization for $\la$-homogenizable polynomials
(see Lemma~\ref{homog}).

(4) Write $\wi f^{(\la)} = (X-b)m$ in $GE[X]$ with $(X-b)\nmid m$.
By Lemma~\ref{homog}, $X-b$ and $m$ are $\la$-homogenizable, 
and since $X-b$ is prime in $GE[X]$, $\text{gcd}(X-b,m) = 1$.  
So by (3) above, $f = gh$  for $\lambda$-polynomials 
$g$ and $h$ in $E[x]$ with $\wi g^{(\la)} = X-b$ and 
$\wi h^{(\la)} =m$.
Write $g = c_1X+c_0$; so, $\wi c_1 = \wi1$ and $\wi c_0 = b$.  
Then, $a = c_0c_1^{-1}$ is a root of $g$, so of $f$,  and $\wi a=
\wi c_0 \wi c_1^{-1} =b$. Since $\wi a$ is not a root of 
$h^{(\la)}$,  $a$ cannot be a root of $h$ by Lemma~\ref{split}.
So, $a$ is a simple root of $f$. 

(5) Since $f'$ is $\lambda$-homogenizable in $GE[X]$, 
 for any $i$ with $a_i\neq 0$,
$v(a_i)+i\lambda=v(a_n)+n\lambda$, hence
$v(a_i)=(n-i)\lambda+v(a_n)$. In particular, we have
$v(a_0)=n\lambda+v(a_n)$, hence
$\lambda=\frac{1}{n}(v(a_0)-v(a_n))$. So, by 
Prop.~\ref{samevalue}(3)
$f$ is a $\lambda$-polynomial. Moreover, we have $\ov
a_i^{(\lambda)}=\widetilde a_i$ for any $i$, so
$\widetilde f^{(\lambda)}=f'$. 
\end{proof}

For monic $\la$-polynomials over a Henselian valued 
field, Lemma~\ref{split} and Th.~\ref{lambda}(2) and (4) were 
essentially proved by Boulagouaz in 
\cite[Lemme~2.4, Th.~2.5, Cor.~2.6]{B(b)98}.

\begin{Cor}\label{lift}
Let $E$ be a Henselian valued field, and let 
$g=\sum_{i=0}^n b_iX^i$ be
 a $\la$-homogenizable irreducible polynomial of $GE[X]$, with
$b_0\ne 0$.  Choose any $a_i\in E$ with $\wi a_i = b_i$, 
$0\le i\le n$, and let $f = \sum_{i=0}^n a_iX^i$.  
Then, for any root $a$ of $f$ in an algebraic extension $K$
of $E$ such that $\wi a$ is a root of $g$ in $GK$, we have
$G(E[a]) = GE[\wi a]$ and $[E[a]:E] = [G(E[a]):GE]$.    
\end{Cor}

\begin{proof}
By Th.~\ref{lambda}(5) $f$ is a $\la$-polynomial in $E[X]$,
with $\wi f^{(\la)} = g$. 
So by Th.~\ref{lambda}(1), $f$ is irreducible in~$E[X]$.  
Now, clearly $GE[\widetilde a]\subseteq G(E[a])$.
But, as $f$ and $g$ are irreducible,  
$$
[GE[\widetilde a]
  : GE] \, = \, 
\deg(g) \, = \, \deg(f) \, = \, [E[a] : E] \, \geq \,  
[G(E[a]) : GE] \, \ge  \,[GE[\wi a]) : GE].
$$ 
Hence, equality holds throughout, which implies that
$G(E[a])=GE[\widetilde a]$.
\end{proof}

\begin{Rq}
 Note that (1) and (5) of Th.~\ref{lambda} are true without 
assuming that $v$ is Henselian.  So, the Henselian assumption 
 can be omitted from  Cor.~\ref{lift}, as well.
\end{Rq}

\begin{Pro} \label{wild}
Let $E$ be a Henselian valued field with residue 
characteristic $p>0$ and $L$ a
 purely wild  finite-dimensional 
graded field extension of $GE$,
 then there is a defectless field extension  $K$
 of $E$ such that $GK=L$. If  ${\rm char}(E)=p$, then  
$K$ can be  chosen to be  purely inseparable field
  extension of $E$.
\end {Pro}

\begin{proof}
 Let $N$ be a  field extension of $E$ such that 
$L\subseteq GN$. Assume first
that $L=GE[\widetilde  a]$ for some $a\in N$,  and let 
$p^n=[L : GE] = [q(L):q(GE)]$. Since $q(L)$ is
 purely inseparable over $q(GE)$,  the minimal polynomial of 
$\widetilde  a$
 over $q(GE)$ is  $g:=X^{p^n}-\widetilde  a^{p^n}$, 
where $\wi a^{p^n} \in 
L^*\cap q(GE)\subseteq GN^*\cap q(GE)=GE^*$. So,  there is $b\in E$ 
such that
   $\wi b = \widetilde  a^{p^n}$. Let $f=X^{p^n}-b$ and 
let  $x$ be a root of $f$ in
   some finite-dimensional field extension $N'$ of $N$. Clearly, 
$\widetilde  x$ is a root of
   $g$; so, $\widetilde  x=\widetilde  a$ because 
$\widetilde  x^{p^n}=\wi b = \widetilde  a^{p^n}$
   in $GN'$.  Note that $g$ is $\gr(\wi a)$-homogenizable, so 
it is irreducible in $GE[X]$ by Lemma~\ref{homog}(3).  
By Cor.~\ref{lift}, we have  
$G(E[x])=GE[\widetilde x]=GE[\widetilde  a]=L$.

Now, let $L$ be an arbitrary finite-dimensional purely wild 
graded field extension  of $GE$.
 Then, we can write 
$L=GE[\widetilde  {a_1}, ..., \widetilde  {a_r}]$, and the result
follows by induction on $r$.  
\end{proof}

As a consequence of Th.~ \ref{lambda},  we have the following
Corollary which gives  explicitly the correspondence between
(finite-dimensional) tame valued field extensions over a 
Henselian
valued field and tame graded field extensions.
Recall that if $L$ is a finite-dimensional extension of a 
Henselian valued field $E$, then $L$ is {\it tame} 
(or tamely ramified) over $E$ if $\text{char}(\ov E) = 0$
or $\text{char}(\ov E) = p >0$, $\ov L$ is separable over 
$\ov E$, $p\nmid |\Gamma_L:\Gamma_E|$, and $[L:E] \,|\Gamma_L:
\Gamma_E| = [L:E]$.

\begin {Cor}\label{corresp} 
\cite[Th.~ 5.2]{HW(a)99} Let $(E, v)$ be a Henselian 
valued field.
Then, the map $K\mapsto GK$ gives a one-to-one correspondence 
between the set of $E$-isomorphism
classes of finite-dimensional tame field extensions of $E$ and 
the set of graded $GE$-isomorphism classes
of finite-dimensional tame graded field extensions of $GE$. 
Moreover, $K$ is a Galois tame
$($finite-dimensional$)$ field extension of $E$ if and only if $GK$ 
is a Galois
$($finite-dimensional$)$ graded field extension of $GE$, in which  
case  $\Gal(K/E)\cong\Gal(GK/GE)$.
\end {Cor}

\begin{proof}
 If $K$ is a tame field extension of $E$, then 
obviously $GK$ is a tame graded
field extension of~$GE$. Let $K'$ be a tame field extension of 
$E$ such that $K'\cong K$.
 Since $E$ is Henselian, the isomorphism
respects the valuations on $K$ and $K'$ extending $v$ on E;
so, ${GK\cong_gGK'}$.

Conversely, if $L$ is a tame finite-dimensional graded field 
extension of $GE$, then $q(GL)$ is separable over 
$q(GE)$ by \cite[Prop.~3.5]{HW(a)99};  so,  we can
 write $L=GE[\widetilde  x_1, ..., \widetilde  x_r]$, where 
$x_i\in E_{\alg}$ with
 $\widetilde  {x_i}$ separable over $q(GE)$. 
Let $g$ be the minimal polynomial of $\wi x_1$ over $q(GE)$.  
Then $g$ is $\la$-homogenizable in $GE[X]$ where 
$\la = \gr(\wi x_1)$,
 $g$ is irreducible in 
$GE[X]$ by Lemma~\ref{homog}(3), and $\wi x_1$ is a simple root
of $g$.  Take any $f = \sum_{i=1}^n c_iX^i \in E[X]$ such that 
$\sum_{i=1}^n \wi c_iX^i = g$.  By Th.~\ref{lambda}(5), $f$ is 
a $\la$-polynomial with $\wi f^{(\la)} = g$. So, 
$f$ is irreducible in  $E[X]$ by Th.~\ref{lambda}(1),  
and Th.~\ref{lambda}(4)
applied over $E[x_1]$ shows that $f$ has a simple root 
$a_1 \in E[x_1]$
with $\wi a_1 = \wi x_1$.  By Cor.~\ref{lift}, $G(E[a_1]) = 
GE[\wi a_1]$.  Moreover, $E[a_1]$ is tame over $E$, as $G(E[a_1])$
is tame over $GE$ and $[E[a_1] :E] = [G(E[a_1]):GE]$ by 
Cor.~\ref{lift}. Since $L = G(E[a_1])[\wi x_2, \ldots, \wi x_r]$,
which is tame over $G\big(E[a_1]\big)$, by induction on $r$ there exist 
$a_2, \ldots, a_r\in E[a_1][x_2, \ldots, x_r]$ such that 
each $\wi a_i = \wi x_i$ and $E[a_1][a_2, \ldots, a_r]$ is 
tame over $E[a_1]$ with $G(E[a_1][ a_2, \ldots, a_r]) =
G(E[a_1])[\wi a_2,\ldots, \wi a_r] = L$.  Let $K = E[a_1, \ldots,
a_r]$.  Then $GK = L$ and $K$ is tame over $E$, as 
$K$ is tame over $E[a_1]$ and $E[a_1]$ is tame over $E$.
For the uniqueness of $K$ up to isomorphism, suppose 
$K^{\prime\prime}$ is another tame field extension of~ 
$E$ such with a graded $GE$-isomorphism $\eta\colon 
L \to GK^{\prime\prime}  $. Let $b = \eta(\wi a_1)$, 
which is a root of the $g$ above in~ 
$GK^{\prime\prime}$.
 With the $f$ above,
Th.~\ref{lambda}(4) applied over $K^{\prime\prime}$ shows that 
$f$ has a root $a_1^{\prime\prime}$ in $K^{\prime\prime}$
with $\wi a_1^{\prime\prime} = b$.
Since $f$ is irreducible in $E[X]$ with roots $a_1$ and 
$a_1^{\prime\prime}$, we have an $E$-isomorphism $\psi\colon
E[a_1] \to E[a_1^{\prime\prime}]$ with $\psi(a_1) = 
a_1^{\prime\prime}$.
The induced $GE$-isomorphism $\wi \psi\colon
G(E[a_1]) \to G(E[a_1^{\prime\prime}])$ maps
$\wi a_1$ to $\wi a_1^{\prime\prime} = \eta(\wi a_1)$.  
 So, $\wi\psi = \eta|_{G(E[a_1])}$,
as $G(E[a_1]) = GE[\wi a_1]$. Since $K = 
E[a_1][a_2, \ldots, a_r]$, it follows by induction 
on $r$ that there is an $E$-isomorphism 
$K \to K^{\prime\prime}$ inducing $\eta$ on the 
graded fields. 

Now, let $K$ be a Galois tame finite-dimensional field extension 
of $E$, let $G = \Gal(K/E)$, and 
let $w$ be the unique extension of $v$ to $K$.
 By  Prop.~\ref{tandn} and Th.~\ref{GMnormal}, $GK$ is a 
Galois 
graded field extension of 
$GE$.   Take any $\sigma\in G$.  Since $w\circ \sigma
=w$, $\sigma$ induces 
a graded $GE$-automorphism 
$\widetilde  \sigma\colon GK \rightarrow GK$ 
satisfying  
   $\wi\sigma (\wi x)= \widetilde{\sigma(x)}$
for all $x\in E$. Let 
$\varphi\colon G\rightarrow \Gal(GK/GE)$ be the group 
homomorphism 
defined by $\varphi(\sigma)=\widetilde \sigma$, and let 
$G^v = \ker(\varphi)$. So, 
$$
G^v \ = \  \{\sigma\in G \,|\ 
\wi{\sigma(x)} = \wi x \text{ for all $x\in K$}\} 
  \  = \  \{\sigma\in G \,|\ w(\sigma(x)- x) > w(x)
\text{ for all $x\in K\setminus \{0\}$}\},
$$
which shows that $G^v$ is the ramification group for $w$
over $E$ (cf.~\cite[Th.~(20.5)(c)]{E72}).  
But, since $K$ is tame over $E$ the ramification 
group is trivial,  e.g., by the table in \cite[p.~171]
{E72} as $K/E$ is defectless; hence, $\varphi$~
is injective.  Since $|\Gal(GK/GE)| = [GK:GE] = [K:E] = |G|$, 
$\varphi$ is a group isomorphism.

Let $M$ be a tame finite-dimensional field extension of $E$ such 
that $GM$ is a Galois graded
field extension of $GE$ and consider a Galois tame 
finite-dimensional field extension $N$ of
 $E$ containing $M$. By the above $GN$ is a Galois graded field 
extension of $GE$ [resp.,
 of $GM$] and $\Gal(N/E)\cong \Gal(GN/GE)$ 
[resp., $\Gal(N/M)\cong \Gal(GN/GM)$]. Since $GM$ is
 a Galois graded field extension of $GE$, then $\Gal(GN/GM)$ is a 
normal subgroup of
  $\Gal(GN/GE)$, therefore $\Gal(N/M)$ is a normal subgroup of 
$\Gal(N/E)$. Hence, $M$ is a
  Galois field extension of $E$.
\end{proof}

Let $L/F$ be a  finite-dimensional Galois graded field extension. 
In the same way as for ungraded fields, one may define the 
norm $N_{L/F}$ of $L$ over $F$ by 
$N_{L/F}=\prod_{\sigma\in \Gal(L/F)}\sigma(x)$ for all $x\in L$. 
The following lemma is the graded version of Hilbert's Th.~90.

\begin {Lem}\label{hilb90} Let $L/F$ be a finite-dimensional  Galois 
graded field extension   with 
cyclic Galois group 
generated by $\sigma$. Then,  for any  $x \in L^*$,  
$N_{L/F}(x)=1 $ 
if and only if there exists $y \in L^*$ such that 
$x=y\sigma(y)^{-1}$
\end {Lem}

\begin{proof}
Since $L$ is Galois over $F$, $q(L)$ is Galois over $q(F)$
with $\Gal(q(L)/q(F)) \cong \Gal(L/F)$.  Hence, for any 
$x\in L$, $N_{L/F}(x) = N_{q(L)/q(F)}(x)$.  
 Assume that $N_{L/F}(x)=1 $. Then, as
$N_{q(L)/q(F)}(x)=1$, by Hilbert's Th.~90
there is $z \in q(L)^*$ such that $x=z\sigma(z)^{-1}$. We may
assume $z \in L\setminus \{0\}$. Write $z=z_1+...+z_r$, where 
 all the $z_i$  are nonzero homogeneous elements of $L$
and  ${\rm gr}(z_i)<{\rm gr}(z_{i+1})$ for all ~$i$, $1 \leq i<r$.
Since $\sigma(z)x=z$ and $x$ is homogeneous,  for every $i$,
 $\sigma(z_i)x=z_i$.  We can take for $y$ any $z_i$. The converse is clear. 
\end{proof}

%% An alternative cohomological proof of Lemma 1.13 can be obtained 
%% by considering the $G$-module $E^*$.

\begin {Pro} Let $F$ be a graded field  and $n$ a positive integer
with ${\rm{char}}(F)\nmid n$. Suppose $F_0$~contains a 
primitive $n$-{th} root of unity $\zeta$. Then, 
\begin{enumerate}
\item[(1)] 
 If $L$ is a cyclic Galois graded field extension of $F$ 
with $[L:F] = n$, then there is
$x \in L^*$ such that $L=F[x]$, $x^n \in F^*$ and $\Gal(L/F)$ is 
generated by the graded
$F$-automorphism $\sigma$ defined by $\sigma(x)=\zeta x$.
\item[(2)] 
Conversely, if $a \in F^*$ and $x$ is a root of the 
polynomial $X^n-a$ in $q(F)_{\alg}$, then $F[x]$ is a cyclic 
Galois graded 
field extension of $F$ with $[F[x]:F]=m$, where  $m\, |\,n$ and  
${x^m \in F^*}$.
\end{enumerate}
\end {Pro}

\begin{proof}
 (1)  Let  $\sigma$ be a generator of $\Gal(L/F)$. 
We have
$N_{L/F}(\zeta ^{-1})=1$, so by Lemma~\ref{hilb90} there is $x \in L^*$ 
such that $\sigma (x)=\zeta x$. Accordingly, 
$\sigma(x^n)=\sigma(x)^n=(\zeta x)^n=x^n$. Hence, $x^n \in F^*$. 
Since $\sigma^i(x)=\zeta^ix$, the $\sigma^i(x)$
are  pairwise 
distinct for $1\le i\le n$, which implies that the minimal 
polynomial of $x$ over $q(F)$ is 
$X^n-x^n$. So, $[F[x]:F] = n =  [L:F]$, showing that $F[x]=L$.

$(2)$  By  \cite[Th.~ 6.2, p.~ 324]{L84}, $q(F)(x)$ is cyclic 
Galois over $q(F)$  and $[q(F)(x): q(F)] = m$ where $m\,|\,n$ and
$x^m \in q(F)^*$. Hence, by \cite[Th.~ 3.11]{HW(a)99} 
$F[x]/F$ is cyclic Galois of dimension $m$, and by
 \cite[Cor.~ 2.5(b)]{HW(a)99},  $x^m \in F[x]^* \cap q(F)=F^*$. 
\end{proof}

 Let $F$ be a graded field with $\text{char}(F) = p>0$. Then, 
Galois graded $p$-extensions of $F$ are inertial over $F$, so 
they are exactly 
graded fields of the form  $KF$, where $K$ is any Galois 
$p$-extensions of $F_0$.
Thus, a graded field extension $L/F$ of dimension a power of $p$ 
is  cyclic if  
$L=F(x_1,..., x_n)$, where $x=(x_1,..., x_n)\in W_n(L_0)$ and 
$(x_1^p,..., x_n^p)-(x_1,..., x_n)\in W_n(F_0)$ (here $ W_n(L_0)$ 
is the ring of  Witt vectors associated to the field $L_0$). In 
particular, cyclic extensions of degree $p$ of $F$ are $F[x]$, 
where $x$ is a root of a polynomial $X^p-X-a$ for some $a\in F_0$ 
with $x\notin F_0$.

\pagestyle{myheadings}
\markright{{\sc Sec.~\ref{Subfields}: Subfields of nondegenerate tame semiramified division algebras}}
\section {Subfields of nondegenerate tame semiramified division algebras }
\label{Subfields}

\setcounter{Pro}{0}

For a central simple algebra $A$ over  a field $E$, as usual we 
set $\deg(A) = \sqrt{[A:E]}$ and $\exp(A) =$ the order of 
$[A]$ in the Brauer group $\text{Br(E)}$.

Before reviewing the notion of nondegeneracy, we recall
Tignol's Dec groups.  Let $N$ be a finite-dimensional Galois field 
extension of a field $E$ with abelian Galois group 
$G = \Gal(N/E)$.  Since $G$ is abelian, there is a base
$(\sigma_1, 
\ldots, \sigma _m)$ of $G$, i.e.,  
$G = \langle \sigma_1 \rangle \oplus \ldots \oplus
\langle \sigma_m\rangle$.  Let $r_i$ be the order of $\sigma_i$
in $G$.  For each~$j$, let $K_j$ be the fixed field of the subgroup 
of $G$ generated by all the $\sigma_i$ for $i\ne j$.  So, $K_j$
is a cyclic Galois extension of $E$ with $[K_j:E] = r_j$ and 
$\Gal(K_j/E) = \langle \sigma_j|_{K_j}\rangle$; also,
$N = K_1\otimes_E \ldots \otimes _E K_m$. The group 
$\text{Dec}(N/E)$, introduced by Tignol in \cite{T},
is the subgroup of the relative 
Brauer group $\text{Br}(N/E)$ (\,$= \ker(\text{Br}(E)
\to\text{Br}(N)$\,) generated by all the subgroups 
$\text{Br}(L/E)$ as $L$ ranges over the fields with 
$E\subseteq L \subseteq N$ and $\Gal(L/E)$ cyclic.  Equivalently, 
$\text{Dec}(N/E)$ is the subgroup of $\text{Br}(E)$ generated
by the $\text{Br}(K_i/E)$ for $1\le i\le  m$.  Tignol
showed in \cite[Cor.~1.4]{T} that $\text{Dec}(N/E)$ consists of the
Brauer classes of central simple $E$ algebras $T$ containing 
$N$ such that $\deg(T) = [N:E]$ and $T$ is a tensor product of 
cyclic algebras with respect to the $K_i$, i.e., 
$T \cong (K_1/E,\sigma_1, c_1) \otimes_E \ldots \otimes_E
(K_m/E, \sigma_m,c_m)$.  Such algebras~$T$ were said by Tignol 
to decompose according to $N$, whence the name $\text{Dec}(N/E)$.
As the definition makes clear, $\text{Dec}(N/E)$ is intrinsic
to $N$ and $E$, and is independent of the choice of cyclic 
decomposition of~$\Gal(N/E)$.   

Now, with $E,N,G$, and the $\sigma_i$ as above, let $A$ be a 
central simple $E$-algebra containing $N$ as a maximal subfield
with $\deg(A) = [N:E]$.  For $1\le i\le m$, choose (by 
Skolem-Noether) $z_i\in A^*$ with $z_icz_i^{-1} = \sigma_i(c)$
for all $c\in N$.  Let $u_{i,j} = z_iz_jz_i^{-1}z_j^{-1}
\in C_A(N)^* = N^*$, and let $b_i = z_i^{r_i}$, which also 
lies in $N$, as $\sigma_i^{r_i} = \text{id}_N$.  Then, 
$A = \bigoplus\limits_{0\le j_1\le r_1-1} \ldots
\bigoplus\limits_{0\le j_m\le r_m-1}Nz_1^{j_1}\ldots z_m^{j_m}$,
and the multiplication in $A$ is completely determined by 
$N$, $G$, the $u_{ij}$, and the $b_i$, so we write 
$A = (N/E, G, S, U, \mathbf b)$, where $U= (u_{i,j})_{1\le i,j\le m}$, 
$\mathbf b = (b_i)_{1\le i\le m}$, and $S = (\sigma_i)_{1\le i\le m}$
is the chosen base of $G$.  Amitsur and Saltman defined in 
\cite[p.~81]{AS78} a condition that they called {\it degeneracy}
for the matrix $(u_{i,j})$, which by \cite[Prop.~0.13]{BM00}
is equivalent to: there is a field $L$, $E \subseteq L\subseteq N$
such that $\Gal(N/L)$ is noncyclic and $[C_A(L)] \in
\text{Dec}(N/L)$.  When there is such an $L$, we say that 
$N$ is {\it degenerate in} $A$, or (when $N$ is understood)
$A$ is {\it degenerate}.  When there is no such $L$, we say 
that $N$ is {\it nondegenerate} in $A$.  Note that the 
characterization in \cite{BM00} makes it clear that 
degeneracy is intrinsic to $N$ and $A$, independent of 
the presentation of $G$ and of the choice of the~ 
$z_i$.
(However, 
degeneracy is not intrinsic to $A$.  Indeed, 
K.~McKinnie has recently given in \cite{Mc3} an example 
of a central division algebra $A$ over a field $E$ with 
maximal subfields $N$ and $N'$ each
abelian Galois over $A$ such that $N$ is degenerate in $A$
but $N'$~is not.)  Clearly, if $\Gal(N/E)$ is cyclic, then
$N$~is nondegenerate in $A$.  Also,
it is easy to see that if $[N:E]$~has more than 
one distinct prime factor, then $N$~is nondegenerate in $A$ if and 
only if each primary component of $N$ is nondegenerate in the 
corresponding primary component of $A$.  Therefore, our focus
will be on nondegenerate algebras of prime power degree
with $\Gal(N/E)$ noncyclic.
The first examples of nondegenerate algebras (with $\Gal(N/E)$
noncyclic) given in \cite[Remark, p.~82]{AS78} satisfied $\exp(A) = \deg(A)$,
(for which the nondegeneracy holds trivially, 
see \cite[Lemma~1.7]{AS78}).  Subsequently, 
Saltman gave in \cite[Cor.~12.15]{S99} an example of a 
nondegenerate 
  algebra $A$ with $\deg(A) = p^2$ and $\exp(A) = p$ for any odd
prime $p$ over a field $E$ containing a primitive $p$-th root of 
unity.  Recently, McKinnie has given in \cite {Mc08} 
more examples of nondegenerate division algebras of 
odd prime exponent, from which further examples can be built
as well. See Remarks~\ref{existence} below.

Now let $F$ be a graded field, and let $B$ be an  
inertially split graded $F$-central division algebra.
Then, as defined in \cite[Remark~2.13]{M08}, $B$ is said to be 
{\it degenerate} if
it has a graded subfield $L$ inertial over $F$ such that $C_B(L)$
is nicely semiramified with $\Gamma_{C_B(L)}/\Gamma_L$ noncyclic.
Assume now that $B$ is semiramified as well as 
inertially split.  Then, $B_0$ is a field abelian Galois over $F_0$
 with 
$\Gal(B_0/F_0)\cong \Gamma_B/\Gamma_F$.  The graded subfields
$L$ of~$B$  inertial over $F$ are  the graded subfields of 
$B_0F$ containing $F$, and are in one-to-one
correspondence with the subfields of $B_0$ containing~ 
$F_0$. (\,$L \leftrightarrow L_0$; note that $L = L_0 F$.)
In particular, $B_0F$ is a maximal graded
subfield of $B$; it is inertial over $F$, and it contains 
all other graded subfields of $B$  inertial over $F$.
Furthermore,  $B_0F$ is Galois over $F$ with   
$\Gal(B_0F/F) \cong \Gal(B_0/F_0)$, which is abelian.
In this context, the degeneracy of $B$ is equivalent to 
what can be called 
the degeneracy of $B_0F$ in~$B$; that is,  $B$ is degenerate iff 
 there is an inertial graded field extension 
$L$ of $F$ in $B$ such that $\Gamma_{C_B(L)}/\Gamma_L$ (\,$\cong 
\Gal(B_0F/L) \cong\Gal(B_0/L_0)$\,) is noncyclic and 
$C_B(L)$ is isomorphic to a tensor product of cyclic graded
algebras over $L$ with respect to cyclic Galois graded 
field extensions of $L$ within $B_0F$. (This equivalence holds 
because $C_L(B)$ is inertially split and semiramified by 
\cite[Prop.~1.3]{M08}, 
so every cyclic subalgebra of $C_L(B)$ determined by 
an inertial cyclic graded field extension of $L$ is nicely   
semiramified, by \cite[Prop.~1.3]{M08}.)
By \cite[Prop.~2.15]{M08},  
$B$ is degenerate if and only if  the maximal subfield $q(B_0F)$
is degenerate in the $q(F)$-central division algebra
$q(B)$. Also, let $\sigma_1, \ldots, \sigma_m$ be any base of 
$\Gal(B_0/F_0)$, and choose any $y_i\in B^*$ with 
$y_icy_i^{-1} = \sigma_i(c)$ for all $c\in B_0$.  Let 
$u_{i,j} = y_iy_jy_i^{-1}y_j^{-1} \in B_0^*$. Then, 
by \cite[Prop.~2.17]{M08}
$B$ is degenerate if and only if the collection
$(u_{i,j})_{1\le i,j\le m}$ satisfies the Amitsur-Saltman 
degeneracy condition as elements of the abelian Galois field
extension $B_0$ of $F_0$.

The generic abelian crossed products of Amitsur and Saltman
are associated to such graded division algebras.  Specifically, 
let $A = (N/E, G, S, U, \mathbf b)$ be an abelian crossed
product  over a field~$E$, as described above.
 From the data associated to $A$, 
Amitsur and Saltman defined in \cite[p.~83]{AS78} a generic  
abelian crossed product
$ A' = \mathcal K(N/E, G, S, U)$ which is a division algebra of the 
same degree as $A$ whose center $Z$  is purely transcendental 
over $E$;  $A'$ has a maximal subfield $M = N\otimes_E Z$
which is abelian Galois over $Z$ with $\Gal(M/Z)
\cong \Gal(N/E) = G$, and $ A^{\prime *}$ contains elements
$y_1, \ldots, y_m$ such that $y_i$ induces $\sigma_i$ on 
$M$ by conjugation, and $y_iy_jy_i^{-1}y_j^{-1} = u_{i,j}$
for all $i,j$.  This $A'$ depends up to isomophism
on the choice of base $S = (\sigma_1, \ldots \sigma_m)$
of $G$ and on $U = (u_{i,j})$, but not on $
\mathbf b = (b_1, \ldots,
b_m)$.  Also, it follows from \cite[Prop.~2.3]{T} 
that $M$ is degenerate in $A'$ iff 
$N$  is degenerate in $A$.  $A'$ is 
definable as the 
ring of quotients of an iterated twisted polynomial ring, 
and it was  shown in \cite[Th.~1.1]{BM00} 
that $A'$ is therefore also $q(B)$ for a 
graded division ring $B$, which is an iterated twisted Laurent
polynomial ring.  Let $F = Z(B)$, a graded field.
It was shown further in \cite[Th.~1.1]{BM00} that  $q(F) =Z$,
 $F_0 = E$, and  
$B$ is inertially split and semiramified over $F$ with 
$B_0 = N$ and $\Gamma_B = \Z^m$. Moreover, by 
\cite[Prop.~2.15]{M08}, $M$ is degenerate in~$A'$ iff
$B$ is degenerate.  We will see in Cor.~\ref{generic}
below how results on subfields of nondegenerate
algebras over Henselian fields  
 yield  another proof of one of Saltman's key results 
about maximal subfields
generic abelian crossed products.

\begin{remarks}\label{existence} 
 (i) Let $p$ be an  odd prime number, and let $G$ be a noncyclic finite 
abelian $p$-group of order~$p^n$, $n\ge 2$.  McKinnie gave in 
\cite[Cor.~3.2.11]{Mc08} an example of a central 
division algebras $A$ over any suitable  field $E$ of any 
characteristic with 
maximal subfield $N$ nondegenerate in $A$ with $\Gal(N/E) \cong G$
and $\exp(A) = p$.  This yields nondegenerate algebras of 
higher exponent, as follows: 
Say $A = (N/E, G, S, U, \mathbf b)$, as above.
Let $A' = \mathcal K(N/E, G, S,U)$ be the associated generic
abelian crossed product. Let $E' = Z(A')$, which is purely 
transcendental over $E$, and let 
$N' = N\otimes _E E'$, which is a maximal subfield which is 
nondegenerate in $A'$ with  $\Gal(N'/E')\cong \Gal(N/E) \cong G$.  
But, 
$\exp(A') = \text{lcm}\big(\exp(A), \exp(G)\big)$ by 
\cite[Prop~2.3]{T}.  Thus, 
simply by choosing $G$ to have exponent $p^r$ for $1\le r\le n-1$, 
we obtain nondegenerate abelian crossed products of exponent 
$p^r$ and degree $p^n$ (cf.~\cite[Ex.~3.3.1]{Mc08}).    

(ii)  From any nondegenerate generic abelian crosed product
$A' = \mathcal K(N/E, G, S,U)$ of degree $p^n$ and exponent
$p^r$, one can obtain a nondegenerate inertially split 
semiramified division algebra $A^{\prime\prime}$ over a
Henselian valued field with $\deg(A^{\prime\prime}) = p^n$
and $\exp(A^{\prime\prime}) = p^r$.  For $A^{\prime\prime}$, 
nondegeneracy is defined to mean nondegeneracy in 
$A^{\prime\prime}$ of the (unique up to isomorphism)
maximal subfield of $A^{\prime\prime}$ which is inertial over~ 
$Z(A^{\prime\prime})$.  One can obtain such an $A^{\prime\prime}$
as what McKinnie calls the \lq\lq power series generic abelian
crossed product \cite[Def.~3.6]{Mc07}, in which the iterated twisted
polynomials in $A'$ are replaced by iterated twisted
Laurent series.  Another way to produce such an $A^{\prime\prime}$
is to view $A'$ as $q(B)$ for $B$ an inertially split graded 
division algebra, and let $A^{\prime\prime} = 
A'\otimes_{Z(A')} HZ(A')$, where $HZ(A')$ is the Henselization
of $Z(A')$ with respect to a valuation on $Z(A')$ induced by 
the grading on $Z(B)$.  (See the proof of Cor.~\ref{generic}
below.) 

(iii) If one wants nondegenerate abelian crossed product 
algebras with specified exponent exceeding~$\exp(G)$, these
are obtainable by a slight adaptation of McKinnie's examples
as follows:  Her $A$~is obtained as 
$\mathcal A\otimes _{\mathcal E}E$, where $\mathcal A$ is a 
division algebra with nondegenerate maximal subfield 
$\mathcal N$ with $\Gal(\mathcal N/\mathcal E) \cong G$
but $\exp(\mathcal A) = \deg(\mathcal A) = |G|$, while
$E$ is the function field $E = \mathcal E(Y)$, where $Y$ is the 
Brauer-Severi variety of~$\mathcal A^{\otimes p}$.  Passage
from $\mathcal E$ to $E$ generically splits 
$\mathcal A^{\otimes p}$, so generically reduces the 
exponent of $\mathcal A$ to $p$ (while assuring that 
$\mathcal A \otimes_{\mathcal E}E$ is a division ring.)  
For any $r$ with $1\le r <n$, let $\mathcal E' = 
\mathcal E(Z)$, where $Z$ is the Brauer-Severi variety
of $\mathcal A^{\otimes p^r}$, and let $\mathcal A' = 
\mathcal A \otimes _{\mathcal E}\mathcal E'$ and 
$\mathcal N' = \mathcal N \otimes_{\mathcal E} \mathcal E'$.
Then,  by Amitsur's theorem \cite[Th.~5.4.1, p.~125]{GS},
$\ker\big(\text{Br}(\mathcal E) \to \text{Br}(\mathcal E')\big)$
is the cyclic group generated by $\big[\mathcal A^{\otimes p^r}\big]$;
so, $\exp(\mathcal A') = p^r$.  Also, $\mathcal A'$ is a 
division algebra (of degree $p^n$) by the Schofield-van den Bergh
index reduction formula \cite[Th.~1.3]{SB} for function fields of
Brauer-Severi varieties.
Clearly, $\mathcal N'$~is a maximal subfield of $\mathcal A'$
which is Galois over~$\mathcal E'$ with $\Gal(\mathcal N'/
\mathcal  E') \cong \Gal(\mathcal  N/\mathcal E) \cong G$. 
Furthermore, $\mathcal N'$~is nondegenerate in $\mathcal A'$.
To see this, let $\mathcal E^{\prime\prime} = \mathcal E'\cdot E$,
 the free 
composite  of $\mathcal E'$ and $E$
over $\mathcal E$; so, $\mathcal E^{\prime\prime}$ is the 
function field over $E$ of the Brauer-Severi variety of 
$A^{\otimes p^r}$.  Since $A^{\otimes p^r}$ is split, 
$\mathcal E^{\prime\prime}$ is purely transcendental over 
$E$, by \cite[Th.~5.1.3, p.~115] {GS}.  Therefore, since  $N$ is nondegenerate
in $A$ (as \text{McKinnie} proved), $N\otimes _E \mathcal E^{\prime\prime}$
is nondegenerate in~$A\otimes _E \mathcal E^{\prime\prime}$, 
which follows from by \cite[Prop.~2.3]{T}.  Then, as $\mathcal A'\otimes_{\mathcal E'}
\mathcal E^{\prime\prime}\cong A\otimes_E 
\mathcal E^{\prime\prime}$ and $\mathcal N'\otimes_{\mathcal E'}
\mathcal E^{\prime\prime}\cong N\otimes_E 
\mathcal E^{\prime\prime}$, $\mathcal N'$ must be 
nondegenerate in $\mathcal A'$. 
\end{remarks}

%%  Let $F$ be a graded field and $D$ a nondegenerate semiramified
%%  graded division algebra of prime power degree $p^n$ over $F$ with
%%  \(\mbox{\rm rk}(\Gal(D_0/F_0))\) $\geq 2$. We recall that it was
%%  shown in \cite[Lemma 3.1]{M08} that if $d\in D^*$ and $d^p\in F$, 
%%  then
%%   $d\in (D_0F)^*$. Moreover, if in addition 
%%  $p=$ \(\mbox{\rm char}(F)\),
%%  then $d\in F^*$. 

\smallskip

Throughout the rest of this section, let 
$E$ be a field with Henselian 
valuation $v$, and 
let  $D$ be a division algebra with center $E$ and 
$[D:E] = p^n$ for some prime number $p$ and some $n\in \N$.  We 
assume further that $D$ is inertially split semiramified with respect to 
the unique extension of $v$
to a valuaton of~$D$.  

There is a distinguished maximal subfield
$N$ of $D$, namely, an inertial maximal subfield.
This $N$ is unique up to isomorphism in $D$, and, since 
$\ov N = \ov D$ is abelian Galois over $\ov E$ 
(see Prop.~\ref{known}(1) below),
$N$~is abelian
Galois over $E$, with $\Gal(N/E) \cong \Gal(\ov D/\ov E)$.
We  assume further that $D$ is nondegenerate, by which is 
meant that $N$ is nondegenerate in~$D$.
Such $D$ and $N$ exist, as we noted in Remark~\ref{existence}(ii).    
The goal of this section is to obtain information 
about subfields of $D$ (containing $E$).   Of course, 
the inertial subfields are known: their isomorphism
classes are in one-to-one correspondence with the 
subfields of~$\ov D$ containing $\ov E$.  The interest, therefore,
is with the noninertial subfields.
We first recall some 
known properties of  $D$  and its subfields
which will be used repeatedly below.

\begin{Pro}\label{known} 
\begin{enumerate}
\item[]
\item[(1)]

$\ov D$ is abelian Galois over $\ov E$ with 
$\Gal(\ov D/\ov E) \cong \Gamma_D/\Gamma_E$.
\item[(2)] 
If  $\Gamma_D/\Gamma_E$ is noncyclic,
 then $D$ has no $($non-trivial$)$ subfield
 totally ramified over~$E$.
\item[(3)]
If
 $K$ is a subfield of $D$ containing $E$ such that
  $\Gal(\overline D/\overline K)$ is noncyclic, then $K$ is 
inertial over~$E$.
\item[(4)]
Let $M$ be a subfield of $D$ with $M$ inertial over $E$, 
and let $C = C_D(M)$.  Then $C$ is inertially split, 
semiramified, and nondegenerate, with $\ov C = \ov D$.
\end{enumerate}
\end{Pro}

\begin{proof}
(1)  
This holds by  \cite[Lemma~5.1]{JW90}, as $D$ is inertially split
with $\ov D$ a field. 
(2) holds by \cite[Prop.~ 3.2]{M08}, and (3) by 
\cite[Prop.~ 3.3]{M08}. (4)  $C$ is inertially split since $D$
is, and by
\cite[Th.~3.1(b)]{JW90} (or by 
embedding $M$ in an inertial lift of $\ov D$ over $E$ in $D$), 
$\ov C = \ov D$. Since $C$ is inertially split with $\ov C$ a field, 
by \cite[Prop.~1.3(1)]{M08} $C$ is semiramified; the nondegeneracy of
$C$ is immediate from the nondegeneracy of $D$.  
\end{proof}

\begin {Th}\label{normalsubfield}
With the hypotheses above, assume further that 
 ${\rm{char}}(\ov E) =
p$  and 
 $\Gamma_D/\Gamma_E$ is noncyclic. Let $K$ be a  subfield of $D$
with $K$ normal over $E$.
 Then, $K$ is  Galois and 
  inertial over $E$. So, there is a subgroup $H$ of
  $\Gal(\overline D/\overline E)$ such that 
${\Gal(K/E)\cong \Gal(\overline D/\overline E)/H}$.
  In particular, if $K$ is a  maximal subfield of $D$
which is Galois over $E$, then
  $\Gal(K/E)\cong \Gal(\overline D/\overline E)$. 
\end {Th}

\begin{proof}
Let $K$ be a  subfield of $D$ which is normal over $E$. Then by 
Th.~\ref{GMnormal},
$GK$ is a normal graded field extension of $GE$. Let 
$L= GK^{\Gal(GK/GE)}$.
By Prop.~\ref{fix}, $L$ is a purely wild graded field
extension of~$GE$.
  Therefore, by \cite[Lemma 3.1]{M08}, $L=GE$. Hence, by 
Prop.~\ref{fix} 
again, $GK$ is a Galois
 graded field extension of $GE$. In particular, 
by \cite[Th.~3.11]{HW(a)99} $GK$ is a tame
 graded field extension of $GE$;  since also  $[GK:GE]$
is  a power of $p$, we must have $\Gamma_{GK} = \Gamma_{GE}$.
Hence,   $GK$ is inertial over $GE$.
 Therefore, $K$ is an inertial valued
 field extension of $E$, and by Cor.~ \ref{corresp}, 
$K$~is a Galois
 field extension of $E$. So, $\Gal(K/E)\cong 
\Gal(\ov K/\ov E)$, which is a homomorphic image of 
$\Gal(\ov D/\ov E)$. The rest is obvious.  
\end{proof}

 \noindent{\bf Remark.} Independently, McKinnie has proved in 
\cite[Th.~ 1.2.1]{Mc08} that if $D$ is a
 semiramified division $p$-algebra over a Henselian valued field 
$E$ with $\overline D$
 separable over $\overline E$ (which is equivalent to saying $D$ 
is an inertially split semiramified
 division algebra over $E$) and with $\overline D/\overline E$ 
not strongly degenerate 
 (See \cite[Def.~ 1.5]{Mc07}), and if $K$ is a Galois subfield of 
$D$, then $\Gal(K/E)$ is of the form
 $\Gal(\overline D/\overline E)/H$, where $H$ is a subgroup of 
$\Gal(\overline D/\overline E)$.

\begin{Cor}\label{cyclic}
  Suppose ${\rm{char}}(\ov E) =p$.  Then
$D$ is a cyclic algebra if and only if $\Gamma_D/\Gamma_E$ is 
cyclic.
\end {Cor}

\begin{proof} 
Recall from Prop.~\ref{known}(1) that  ${\Gamma_D/\Gamma_E\cong
\Gal(\ov D/\ov E)}$.  Thus, if $D$ is a 
cyclic algebra, 
then Th.~\ref{normalsubfield} shows that $\Gal(\ov D/\ov E)$
is a cyclic group, so $\Gamma_D/\Gamma_E$ is also cyclic.
Conversely,  if $\Gamma_D/\Gamma_E$ is cyclic, 
then \lq the' inertial lift of  
  $\overline D$ over $E$ in $D$ 
 (see \cite[Th.~2.9]{JW90}) is a cyclic maximal subfield of
  $D$. 
\end{proof}

\begin {Cor}\label{generic}  \cite[Th.~ 3.2]{S78} 
Let $\mathcal E$ be a field with ${\rm char}(\mathcal E)
= p\ne 0$, let $\mathcal N$ be a noncylic abelian Galois
field extension of $\mathcal E$ with $[\mathcal N:\mathcal E]=
p^n$, $n \ge 2$, and let $G = \Gal(\mathcal N/\mathcal E)$.
Let $\mathcal  A = \mathcal K(\mathcal N/\mathcal E, G,S,U)$ be any
associated generic abelian crossed product
with $U= (u_{i,j})$ nondegenerate; let 
$Z = Z(\mathcal A)$.  For any subfield $L$ of $\mathcal A$ with
$L$  Galois over $Z$, there is a subgroup $H$ of
$G$ such that $\Gal(L/Z)\cong G/H$. In particular, if $L$ is a
Galois maximal subfield of $\mathcal  A$, then $\Gal(L/Z)\cong G$.
\end {Cor}

\begin{proof}
As recalled preceding  Remark~\ref{existence} above, by 
\cite[Th.~1.1]{BM00}, $\mathcal  A = q(B)$ for some
 semiramified  graded division algebra $B$ with 
$B_0 = \mathcal  N$, $Z(B)_0 = \mathcal E$, and $q(Z(B)) =Z$.
Let $HZ$ be the Henselization of~$Z$ with respect to the 
canonical valuation on $Z$ determined by the grading on $Z(B)$
with some chosen total ordering of $\Gamma_{Z(B)}$; let 
$H\mathcal A = \mathcal A\otimes _Z HZ$.  Because the valuation on 
$Z$ extends to $\mathcal A = q(B)$, by Morandi's Henselization 
Theorem \cite[Th.~ 2]{Mor89} $H\mathcal A$ is a 
division ring with valuation extending the valuation on $HZ$,
and $\ov{H\mathcal  A} \cong \ov{\mathcal A}\cong B_0\cong 
\mathcal N$ and $\Gamma_{H\mathcal A} = \Gamma_{\mathcal A}$;
also, $\ov {HZ}\cong \ov  Z \cong Z(B)_0\cong 
\mathcal E$.  So, $H\mathcal A$ is inertially split
with $\Gal(\ov{H\mathcal A}/\ov {HZ})\cong \Gal(\mathcal N/
\mathcal E) = G$ and associated graded ring $GH\mathcal A
\cong G\mathcal A \cong B$.  Since $U = (u_{i,j})$ is 
nondegenerate in $\mathcal N$ for $G$ with base $S$, by 
\cite[Prop.~2.17]{M08}  $B$ is a nondegenerate graded division
algebra, so by \cite[Lemma~2.14]{M08} $H\mathcal A$ is 
nondegenerate.  Let $L$ be a subfield of $\mathcal A$ with $L$ 
Galois over $Z$.  Then $L\otimes_Z HZ$ is a subfield 
of $H\mathcal A$  Galois over $HZ$ with 
$\Gal((L\otimes_Z HZ)/HZ)\cong \Gal(L/Z)$. 
So, by 
Th.~ \ref{normalsubfield}, 
there is a subgroup $H$ of
 $G$ such that $\Gal(L/Z)\cong G/H$. In particular, if $L$ is a
 Galois maximal subfield of $\mathcal A$, then 
$\Gal(L/Z)\cong G$.
\end{proof}

\begin {Pro}\label{elabelian}  
 Assume  
$\Gamma_D/\Gamma_E$ is noncylic,
and let $K$ be a subfield of $D$
which is elementary abelian Galois over $E$. Then, $K$ 
is inertial over $E$.
Therefore, $D$ is an elementary abelian crossed product if and 
only if
$\Gal(\overline D/\overline E)$ is elementary abelian.
\end  {Pro}

\begin{proof}
  Since $K$ is an elementary abelian field extension 
of $E$,  we can write
 $K=K_1\otimes_EK_2\otimes_E...\otimes_EK_r$, where each $K_i$   
is a cyclic  field extension of
 $E$ with $[K_i : E] = p$.  
By Prop.~\ref{known}(2), $D$ 
contains no proper totally ramified field extensions of 
$E$. Hence, each $K_i$ is inertial over $E$;  so, $K$ is 
also inertial over 
$E$. If $K$ is in addition a
  maximal subfield of $D$, then $\overline D \ (=\overline K)$ is 
elementary abelian over
  $\overline E$. Conversely, suppose that 
$\Gal(\overline D/\overline E)$ is elementary
  abelian and let $M$ be the inertial lift of $\overline D$ over 
$E$ in $D$. Then, $M$~is a
    maximal subfield of $D$ which is Galois over 
$E$ with $\Gal(M/E) \cong\Gal(\ov M/\ov E)$. 
\end{proof}

\begin{Rq}\label{cycliccase}
 Prop.~ \ref{elabelian} is not true if
$\Gamma_D/\Gamma_E$ is cyclic.  
Indeed, let $k$ be a field containing a primitive
$p$-th root of unity, with a cyclic Galois extension 
$L$ with $[L:k] = p^2$. Let $N$ be the field with 
$k\subsetneqq N\subsetneqq L$.  Let $E$ be the Laurent series 
field  $k((X))$, and let $D$ be the cyclic algebra 
$\big(L((X))/E, \sigma, X\big)$, where $\sigma$ is any generator
of $\Gal(L((X))/E)$.  Then, $D$ is a division algebra of degree
$p^2$ with center $E$.  
Moreover, with respect to the Henselian $X$-adic valuation on 
$E$, $D$ is tame and semiramified, with $\ov D = L$.  
This $D$~is 
trivially nondegenerate since $\Gal(\ov D/\ov E)$ is cyclic.
We have $\Gamma_D/\Gamma_E\cong \Gal(\ov D/\ov E) \cong
\Z\big/ p^2\Z$.  But, take $t$ in $D$ with $tat^{-1} = 
\sigma(a)$ for all $a\in L((X))$ and $t^{p^2} = X$.
Then, $D$ also contains the maximal subfield $N((X))[t^p]$ 
which is elementary abelian Galois over $E$.
\end{Rq}

For a finite abelian $p$-group $P$, $\text{rk}(P)$ denotes 
the number of summands in a cyclic decomposition of~ 
$P$; so $\text{rk}(P) = \dim_{\Z/p\Z}(P/pP)$.

\begin{Pro}\label{rk3}
  Suppose
${\rm rk}(\Gamma_D/\Gamma_E)\geq 3$, 
 and let $K$ be a  subfield of $D$ which is abelian 
Galois over~$E$. Then,  $K$ is inertial 
over $E$.
\end {Pro}

\begin{proof}
 Write $K=K_1\otimes_EK_2\otimes_E\ldots\otimes_EK_r$,
where each $K_i$  is a cyclic Galois field extension of $E$.
  So, $\overline {K_i}$~is cyclic over
$\overline E$. Therefore, $\overline D$ cannot be cyclic over
$\overline {K_i}$ (since ${\rm rk}(\Gal(\overline D/\overline
E))\geq 3$). So, by Prop.~\ref{known}(3) each $K_i$ is inertial over
$E$. Hence, $K$ is inertial over $E$. 
\end{proof} 

\begin {Lem}\label{primeexp}
Let $K$ be any subfield of $D$ containing $E$ with 
$K$ not inertial over $E$. Then, 
$\Gal(\ov D/\ov K)$~is cyclic and 
$\Gamma_K/\Gamma_E$ is isomorphic to a subgroup of 
$\Gal(\ov D/\ov K)$.   
In particular,
 if~${\rm exp}(\Gamma_D/\Gamma_E) = p$, then $K$ is a 
maximal subfield of $D$.
\end {Lem}

\begin{proof} 
  Let $M$ be the maximal unramified  
extension of $E$ in $K$, and let $C = C_D(M)$.
By Prop.~\ref{known}(4) above, $C$ is a nondegenerate
inertially split semiramified division algebra with 
$\ov C = \ov D$.
Moreover 
   $K$ is a
subfield of $C$ which is a
non-trivial totally ramified extension of
 the center $M$ of $C$. So, by
Prop.~\ref{known}(1) and (2) applied to  $C$ as an $M$-algebra, 
 $\Gal(\overline C/\overline M)$ is cyclic. 
Then, $\Gal(\ov D/\ov K)$ is cyclic, as
 $\ov D = \ov C$ and $\ov K = \ov M$. 
The 
canonical isomorphism $\Gamma_D/\Gamma_E \cong
\Gal(\ov D/\ov E)$ of Prop.~\ref{known}(1)
is induced by conjugation, 
so it injects $\Gamma_K/\Gamma_E$ into 
$\Gal(\ov D/\ov K)$.  
If  
$p = {\rm exp}(\Gamma_D/\Gamma_E)=
{\rm exp}(\Gal(\overline D/\overline E))$,
 then  the cyclic group $\Gal(\ov C/\ov M)$
has exponent and hence order $p$.
Then,   $p = |\Gal(\ov C/\ov M)| =\text{deg}(C)$,
as $C$ is semiramified.
So the proper extension $K$ of~$M$ is a maximal subfield
of $C$, and hence $K$ is also a maximal subfield of~$D$. 
\end{proof}

Let $H$ be a finite nonabelian group. We say that $H$ is a 
{\it quaternion
group} if $H$ has order $8$ and is generated by two elements $a$
and $b$  satisfying the  conditions  $a^4=b^4=1$,
$a^2=b^2$ and $ba=a^{-1}b$. If $K/E$ is a normal [resp., Galois]
field extension with a quaternion Galois group, we say that 
$K$ is
a quaternion normal [resp., Galois] field extension of $E$. We say
that a finite group $H$ is {\it Hamiltonian} if $H$ is  
nonabelian and 
every
subgroup of $H$ is normal. 
Recall that  a Hamiltonian group is
the direct product of a quaternion group with an abelian group 
of odd order and an abelian group of
exponent two \cite[Th.~ 12.5.4]{Ha59}.

\begin {Th} \label{pexp} 
 Suppose
${\rm exp}(\Gamma_D/\Gamma_E)=p$, and let  
$K$ be a subfield of $D$ containing $E$. Then, 
\begin{enumerate}
\item[(1)] 
If $K$ is not inertial over $E$, then 
 $K$ is a maximal subfield  $D$.
\item[(2)] 
If $\Gamma_D/\Gamma_E$ is noncyclic and  $K$ is a   non-quaternion 
normal maximal subfield of $D$, then either $K$~is  cyclic 
Galois over $E$ with 
 $[K:E] = p^2$ or $K$ is inertial and elementary abelian Galois
over~$E$.
\end{enumerate}
%% Moreover,  in  case  $(1)$ and in the second part of case  $(2)$, 
%% $K$ is elementary abelian Galois over $E$.
\end {Th}

\begin{proof}
  (1) This follows by Lemma~\ref{primeexp}.

(2) Let $G = \Gal(K/E)$ and let 
$I= K^G$. Then, 
$I$ is purely inseparable over $E$, so $\ov I$
is purely inseparable over $\ov E$.  
Since $\ov D$ is separable over $\ov E$, 
we have 
$\overline I=\overline E$; hence $I$ is totally ramified over $E$.
Since $D$ is nondegenerate, Prop.~\ref{known}(2) shows that 
    $I=E$. Therefore, $K$ is Galois over $E$.
 If $L$ is any proper subfield of $K$, then by part (1) of this 
proof $L$ is inertial over $E$,
  hence $L$ is an abelian Galois field extension of $E$. Therefore, any 
subgroup of $G$ is
  a normal subgroup. So, $G$ is Hamiltonian or abelian.
If $K$ is inertial over $E$, then $K$ is Galois over $E$ since
$\ov K$ is Galois over $\ov E$, and $\Gal(K/E)\cong \Gal(\ov K/
\ov E)$; this group is elementary abelian since $\Gal(\ov D/
\ov E)\cong \Gamma_D/\Gamma_E$, which is elementary abelian.
Suppose now $K$ is not inertial over~$E$, and let 
$M$ be the maximal unramified extension of $E$ in $K$.  
We have seen that every proper subfield of~$K$ is inertial
over $E$, so lies in $M$.  Therefore, $\Gal(K/M)$ is the unique 
minimal proper subgroup of~$G$.  So, $G$ admits no nontrivial
direct product decompositions, and since $G$~is assumed 
non-quaternion, it cannot be Hamiltonian.  Hence, $G$ is 
abelian, and since it has a unique minimal proper subgroup
it must be cyclic. 
The exponent assumption then
  implies that the cyclic groups $\Gamma_K/\Gamma_E$ and 
$\Gal(\overline K/\overline E)$ have
  order at most $p$.  So $[K:E]$ is 
at most $p^2$. The fact that
   $K$ is a maximal subfield of $D$ and $\Gamma_D/\Gamma_E$ is 
noncyclic imply that
    $[K : E]=p^2$.
%%    The rest of the proof is obvious since 
%% $\Gal(\overline D/\overline E)\cong \Gamma_D/\Gamma _E$, 
%% which  is elementary abelian.
\end{proof}
 
\begin{Lem}\label{normalval}
Let $K$ be a tame finite-dimensional Galois extension of the
 Henselian field $E$, and let $L$~be a field with  
$E\subseteq L\subseteq K$,
such that $\ov L = \ov K$.  Then, $L$ is Galois over $E$.
\end{Lem}

\begin{proof}
Let $M$ be the maximal unramified extension of $E$ in $L$.
Because $\ov L = \ov K$, $M$ is the maximal unramified extension
of $E$ in $K$.  Since $K$ is Galois, totally ramified, and tame
over $M$, we have the canonical isormorphism
$\Gal(K/M)\cong \text{Hom}(\Gamma_K/\Gamma_M,\Omega)$, 
where $\Omega$ is the group of roots of unity in $\ov K$,
by \cite[(20.11)]{E72}.  For any field $N$ with 
$M\subseteq N\subseteq K$, this isomorphism maps
$\Gal(K/N)$ to $\text{Hom}(\Gamma_K/\Gamma_N,\Omega)$.
This yields a one-to-one correspondence between subgroups
of $\Gal(K/M)$ and subgroups of $\Gamma_K/\Gamma_M$.  Therefore,
for any $\sigma \in \Gal(K/E)$, since $\sigma(L) \supseteq
\sigma(M) = M$ and 
$\Gamma_{\sigma(L)}=
\Gamma_L$, we must have $\sigma(L) = L$. As $K$ is Galois
over $E$, it follows that $L$ is normal, hence also Galois, over $E$.   
\end{proof}

\begin{Def}\label{maxcyclicdef} 
Let $G$ be an abelian group and $H$ a
non-trivial cyclic subgroup of $G$. We say that $H$~is 
{\it maximally
cyclic} in $G$ if there is no cyclic subgroup $H'$ of $G$ such 
that
$H\subsetneqq H'$. 
\end{Def}

Suppose
$\Gamma_D/\Gamma_E$ is noncyclic. We 
 have previously seen in the proof of 
Th.~ \ref{normalsubfield} and 
also in the proof of Th.~ \ref{pexp} 
  that if $K$ is a normal field extension of $E$ in $D$, then $K$ 
is a Galois field extension of $E$.
   From the nondegeneracy of $D$ over $F$
we have seen in Prop.~\ref{known}(3)
that if 
$\Gal(\overline D/\overline
   K)$ is noncyclic, then $K$ is inertial over $E$. In the next
   proposition we will study the case where $\Gal(\overline
   D/\overline K)$ is maximally cyclic in~$\Gal(\overline
   D/\overline E)$.

\begin{Pro} \label{maxcyclic}
 Let
  $K$  be a subfield of
$D$ which is normal and tame over $E$ 
such that $\Gal(\overline D/\overline K)$ is cyclic. Suppose
$\Gal(\overline D/\overline K)$ is maximally cyclic in
$\Gal(\overline D/\overline E)$. Then, $K$ is Galois over $E$.
Furthermore,
\begin{enumerate}
\item[(1)] 
If ${\rm deg}(D)$ is odd, then $K$ is an abelian  field 
extension of $E$.
\item[(2)] 
 If $\Gamma_D/\Gamma_E$ is noncyclic and ${\rm deg}(D)$ is a
power of $2$, then $K$ is either a quaternion or an abelian field
extension of $E$.
\end{enumerate}
Therefore, if ${\rm rk}(\Gamma_D/\Gamma_E)\geq 3$ and $K$ is not 
a quaternion field 
extension of $E$, then $K$ is inertial over $E$.
\end {Pro}

\begin{proof} 
Let $G = \Gal(K/E)$.
By the same argument as in the proof  of 
Th.~ \ref{pexp}(2), 
 $K$ is  Galois over  $E$.
Let $L$ be a field with $E\subseteq L\subseteq K$.  If 
$\ov L=\ov K$, then by Lemma~\ref{normalval} $L$ is Galois 
over $E$.  On the other hand, if $\ov L \subsetneqq \ov K$
then $\Gal(\ov D/\ov L)$ is not cyclic by the maximal cyclicity
assumption.  Therefore, by Prop.~\ref{known}(3),
 $L$~ is inertial, and hence Galois, over $E$.  Since
$\Gal(L/E)$ is Galois in all cases, $G$ is Hamiltonian or abelian.
In case (1), since $|G|$ is odd, $G$ must be abelian. 

(2)  Suppose $G$ is Hamiltonian.  Since $G$ is nonabelian, 
$K$ is not inertial over $E$.  Suppose we have a tensor decomposition 
$K = L_1\otimes_E L_2$ with fields $L_i$ containing $E$.  
Since $K$ is not inertial over $E$, one of the $L_i$, say $L_1$,
is not inertial over $E$.  Then, $\ov {L_1}$ = $\ov K$, as we saw 
above.
Let $M_i$ be the maximal unramified extension of $E$ in $L_i$.
Then, $\ov{M_2} = \ov{L_2} \subseteq \ov K = \ov {L_1} = 
\ov {M_1}$.  Therefore, $M_2$ is  a subfield of $M_1$.
Since $M_1\otimes_E M_2$ is a field (a subfield of 
$L_1\otimes_E L_2$), we must
have $M_2 = E$.  Therefore, $L_2$ is totally ramified over $E$.
Hence,  $L_2= E$, by Prop.~\ref{known}(2).  This shows 
that $K$ admits no nontrivial tensor product decompositions over 
$E$, and hence that $G$ admits no nontrivial direct product 
decompositions.  Therefore, the Hamiltonian group $G$ must be quaternionic, 
so $K$ is quaternionic over~$E$.

The rest of the proposition follows by Prop.~ \ref{rk3}.
\end{proof}

\begin {Cor}\label{odd} 
 Assume ${\rm
deg}(D)$ is odd and ${\rm char}(\overline E)\nmid
{\rm deg}(D)$. Let  $K$ be  a  maximal subfield of $D$
such that $K$ is Galois over $E$.
If $|\Gamma_K : \Gamma_E|={\rm exp}(\Gamma_D/\Gamma_E)$, then 
$\Gal(K/E)$ is abelian and ${{\rm rk}(\Gamma_D/\Gamma_E)\leq 2}$.
\end {Cor}

\begin{proof} Assuming $D\ne E$, the assumption on $\Gamma_K$
assures that $K$ is not inertial over $E$.
Therefore, by Lemma~\ref{primeexp}, $\Gal(\ov D/\ov K)$ is 
cyclic.  Let $M$ be the maximal unramified extension 
of $E$ in $K$, and let $C = C_D(M)$.  By Prop.~\ref{known}(4),
$C$ is inertially split and semiramified with $\ov C = \ov D$.  
Now, $K$ is a maximal subfield of $C$, since it is maximal in $D$, 
and $K$ is totally ramified over $M$.  Since 
$$
|\Gamma_K:\Gamma_M|
 \ = \  [K:M]  \ = \  \text{ind}(C)  \ = \  |\Gamma_C:\Gamma_M|, 
$$ 
we have $\Gamma_K
=\Gamma_C$.  As $\ov M = \ov K$ and $\Gamma_M = \Gamma_E$, 
Prop.~\ref{known}(1) applied to 
$C$ shows that 
$$
\Gal(\ov D/\ov K)  \ = \  \Gal(\ov C/\ov M)
 \ \cong \  \Gamma_C/\Gamma_M  \ = \  \Gamma_K/\Gamma_E.
$$
Therefore, 
  $$
\big|\Gal(\overline D/\overline
K)\big| \ = \ |\Gamma_K : \Gamma_E| \ = \  {\rm exp}(\Gamma_D/\Gamma_E) 
 \ = \  {\rm exp}\big(\Gal(\overline D/\overline
E)\big);
$$
 so, $\Gal(\overline D/\overline K)$ is maximally cyclic in
$\Gal(\overline D/\overline E)$. Hence, by 
(1) and the last assertion of Prop.~ \ref{maxcyclic},  $K$~is
 abelian Galois over $E$ and  ${\rm rk}(\Gamma_D/\Gamma_E)\leq 2$ 
(as $K$ is not inertial  over~$E$). 
\end{proof}

\begin{Rq} Suppose  ${\rm{char}}(\overline E)
\ne p$.   Let $x$ be an element of $D^*$ such that
${{\rm{ord}}(v(x)+\Gamma_E)={\rm{exp}}(\Gamma_D/\Gamma_E)}$, and 
let $L$ be a maximal
subfield of $D$ containing $x$.  Then by Cor.~ \ref{odd},  $L$
cannot be Galois over $E$ if rk$(\Gamma_D/\Gamma_E)\geq 3$.
\end{Rq}

\begin {Pro}\label{inertial} 
Suppose  ${\rm char}(\overline E) 
\ne p$
and  ${\rm rk}(\Gamma_D/\Gamma_F)\geq 3$.  
 Suppose $K$ is a 
subfield of $D$  which is Galois 
but not inertial over $E$.  Then,
$[K : E]\ge p\, {\rm deg}(D)\,{\rm exp}(\Gamma_D/\Gamma_E)^{-1}$.
\end {Pro}

\begin{proof} 
We may assume that $K$ is minimal in $D$ with the property that
$K$ is Galois but not inertial over~$E$.  Let $M$ be the maximal
unramified extension of $E$ in $K$.  Then, $K$ is Galois,
totally ramified,  
and tame over $M$ (as $\text{char}(\ov E) \ne p$).
Let $L$ be a field with 
 $M\subseteq L \subsetneqq K$.
Then, 
$\ov L = \ov K$ since $\ov M = \ov K$,
%% $\ov K = \ov M \subseteq \ov L\subseteq \ov K$, 
so $L$ is Galois over $E$  by Lemma~\ref{normalval}.
Hence,  $L$ is inertial over $E$ by the minimality of~ 
$K$, i.e., $L = M$.  Thus, $M$ is a maximal proper subfield
of~$K$.  Since $K$ is Galois over $M$, this 
 implies that $[K:M] = p$. So, 
$|\Gamma_K:\Gamma_E| = |\Gamma_K:\Gamma_M| = [K:M] = p$.  
Now, by Prop.~\ref{known}(3), $\Gal(\ov D/\ov K)$ is a cyclic 
group. Hence, $[\ov D:\ov K] \le \text{exp}(\Gal(\ov D/\ov E))
= \text{exp}(\Gamma_D/\Gamma_E)$ (see Prop.~\ref{known}(1)).
Therefore, 
$$
[K:E]  \ = \  [\ov K:\ov E]\, |\Gamma_K:\Gamma_E| \ = 
 \ [\ov D:\ov E]\,[\ov D:\ov K]^{-1}p  \ \ge  \ p \, \text{deg}(D) \,
\text{exp}(\Gamma_D/\Gamma_E)^{-1}.
$$
\end{proof}

\end {document}